\documentclass[11pt,a4paper]{article}
\pdfoutput=1
\usepackage[utf8]{inputenc}
\usepackage{amsmath,amsthm,amssymb,amsfonts,mathtools}
\usepackage{authblk}
\usepackage{url,xcolor,bm,enumerate}
\usepackage{graphicx}
\usepackage[margin=1.475in]{geometry}
\usepackage{pgfplots}
\pgfplotsset{compat=1.14}
\usepgfplotslibrary{groupplots}
\usetikzlibrary{arrows.meta,decorations.markings}
\usepackage{caption}
\usepackage{subcaption}
\newtheorem{theorem}{Theorem}
\newtheorem{proposition}[theorem]{Proposition}
\newtheorem{lemma}[theorem]{Lemma}

\theoremstyle{definition}

\newtheorem{remark}{Remark}

\newcommand{\N}{\mathbb{N}}
\newcommand{\Z}{\mathbb{Z}}

\newcommand{\R}{\mathbb{R}}
\newcommand{\C}{\mathbb{C}}

\DeclareMathOperator{\linspan}{span}
\DeclareMathOperator{\re}{Re}
\DeclareMathOperator{\im}{Im}
\DeclareMathOperator{\supp}{supp}

\newcommand{\J}{\mathcal{J}}
\newcommand{\uhat}{\widehat{u}}
\newcommand{\what}{\widehat{w}}

\newcommand{\rhat}{\widehat{r}}
\newcommand{\xhat}{\widehat{x}}
\newcommand{\inv}{^{-1}}
\newcommand{\overbar}[1]{\mkern 1.5mu\overline{\mkern-1.5mu#1\mkern-1.5mu}\mkern 1.5mu}

\newcommand{\Bbar}{\overbar{B}_R}
\newcommand{\dbar}{\overbar{\partial}}
\newcommand{\zbar}{\overline{z}}
\newcommand{\stilde}{\widetilde{\rho}}
\newcommand{\Ssigma}{\Sigma(\rho)}

\newcommand{\M}{\mathcal{M}}
\newcommand{\B}{\mathcal{B}}
\newcommand{\PP}{\mathcal{P}}

\title{Analysis of a Dynamical System Modeling Lasers and Applications for Optical Neural Networks}
\author[1]{Lauri Ylinen\thanks{lauri.ylinen@helsinki.fi}}
\author[2]{Tuomo von Lerber\thanks{t.vonlerber@skoltech.ru}}
\author[2]{Franko K\"{u}ppers\thanks{f.kueppers@skoltech.ru}}
\author[1]{Matti Lassas\thanks{matti.lassas@helsinki.fi}}
\affil[1]{Department of Mathematics and Statistics, University of Helsinki, Finland}
\affil[2]{Center for Photonics and Quantum Materials, Skolkovo Institute of Science and Technology (Skoltech), Moscow, Russian Federation}

\begin{document}
\maketitle

\begin{abstract}
  \noindent
  An analytical study of dynamical properties of a semiconductor laser
  with optical injection of arbitrary polarization is presented. It is
  shown that if the injected field is sufficiently weak, then the
  laser has nine equilibrium points, however, only one of them is
  stable. Even if the injected field is linearly polarized, six of the
  equilibrium points have a state of polarization that is
  elliptical. Dependence of the equilibrium points on the injected
  field is described, and it is shown that as the intensity of the
  injected field increases, the number of equilibrium points
  decreases, with only a single equilibrium point remaining for strong
  enough injected fields. As an application, a complex-valued optical
  neural network with working principle based on injection locking is
  proposed.

  \text{\bf{Keywords}}: dynamical system, semiconductor laser, laser
  with optical injection, complex-valued neural
  network, equilibrium point, stability, bifurcation analysis\\\\
  \text{MSC2020}:
  37N20, 
  34C15, 
  78A60 
\end{abstract}

\section{Introduction}

Self-sustained oscillatory systems will synchronize with an external
source of periodic perturbation, given that the frequency and the
strength of the injection occur within the locking range. A laser
subject to external optical injection behaves the
same~\cite{lau2009enhanced}. What sets optical oscillators apart from
the electronic ones is the nature of propagating electromagnetic field
that has two orthogonal polarization modes which can be observed with
a pair of base polarization components (meaningful reference
coordinate system), be it linear, circular, or some elliptical. In
following treatment, we choose to express polarization in terms of a
complex amplitude $E=(E_-, E_+)\in\C^2$ that multiplies carrier wave
of the form $e^{-i ( k x - \omega t )}$, where $k$ is the wave vector,
$x$ is the spatial coordinate, $\omega$ is the angular frequency, and
$t$ is the time, such that ${k,x,\omega,t}\in\R$. Coordinates $E_\pm$
of $E$ are the \emph{right} $(+)$ and \emph{left} $(-)$
\emph{circularly polarized} components, they are related to the
orthogonal linear components $E_x$ and $E_y$ of the electric field by
\begin{equation*}
  E_x = \frac{E_++E_-}{\sqrt{2}}\text{ and }
  E_y = -i\frac{E_+-E_-}{\sqrt{2}}.
\end{equation*}
Electric field emitted by a laser is
\begin{equation*}
  \mathcal E(x,t)=\re\left( E(t)e^{-i(kx-\omega t)} \right), 
\end{equation*}
where $E(t)$ is called a \emph{slowly varying amplitude}.

In absence of laser cavity anisotropies, the temporal behavior of a
semiconductor laser under external optical injection can be expressed
with a spin-flip rate
equations~\cite{san1995light,martin1997polarization} that describe the
complex-valued components $E_\pm(t)$ of the slowly varying amplitude
$E(t)$ as
\begin{subequations}
  \label{eq:Martin-Regalado}
  \begin{align}
    \frac{d}{dt} E_\pm(t)
    &=\kappa(1+i\alpha)\, (N(t) \pm n(t)-1) E_\pm(t) + \kappa \eta u_\pm(t),\\
    \frac{d}{dt} N(t)
    &= -\gamma(N(t)-\mu) - \gamma(N(t)+n(t))|E_+(t)|^2
      - \gamma(N(t)-n(t))|E_-(t)|^2 ,\\
    \frac{d}{dt} n(t)
    &= -\gamma_sn(t) - \gamma(N(t)+n(t))|E_+(t)|^2
      + \gamma(N(t)-n(t))|E_-(t)|^2 ,
  \end{align}
\end{subequations}
where $N(t)$ and $n(t)$ are real-valued functions; $N$ is the
difference between the normalized upper and lower state populations,
i.e., the normalized total carrier number in excess of its value at
transparency; $n$ is the normalized imbalance between the population
inversions (in reference to the populations of the magnetic
sublevels), $u_\pm$ are the circularly polarized components of the
electric field of an external injection $u=(u_-, u_+)\in\C^2$, that
is, the amplitude of the external light that goes into the laser,
$\eta$ is the coupling efficiency factor, $\alpha$ is the linewidth
enhancement factor that refers to saturable dispersion (Henry factor),
$\mu$ is the normalized injection current, $\kappa$ is the decay rate
of the cavity electric \emph{field} whence $(2\kappa)^{-1}$ is the
cavity photon lifetime, $\gamma$ is the decay rate of the total
carrier number, and $\gamma_s$ is the excess in the decay rate that
accounts for the mixing in the carriers with opposite spins.

The rate equations~\eqref{eq:Martin-Regalado} are derived to model and
explore polarization properties of Vertical-Cavity Surface-Emitting
Lasers (VCSELs). The rate equations use a normalized injection current
such that the unitless injection $\mu$ of $1$ refers to the laser
threshold operation, and $\mu\approx 3$ refers to the output emission of
$1$ mW on a typical VCSEL. In the physical world, an array of VCSELs
is produced on a semiconductor wafer, where stacks of dielectric
materials form high-reflectivity Bragg mirrors on the top and bottom
sides of the wafer. The mirrors confine an active region in between,
comprising just a few quantum wells with a thickness of some tens of
nanometers. Depending on the active region diameter, the threshold
current and the maximum emission power may be tailored for specific
applications.

Lasers are known to exhibit a rich dynamical behavior under external
optical injection~\cite{wieczorek2005dynamical, kane2005unlocking,
  MR2346863, MR2313730, al2013dynamics}.  Depending on laser
properties and the injected optical power and its frequency, the
differential equation system may converge toward an \emph{equilibrium
  point} (a time independent solution, also called \emph{steady
  state}, \emph{stationary point}, or \emph{critical point}) with
locked phase synchronization. This phenomenon is called
\emph{injection locking}~\cite{siegman1986lasers}. Alternatively, the
system may manifest periodic oscillations, or
chaos~\cite{thornburg1997chaos, MR2299634}. In this work, we explore
equilibrium points of system~\eqref{eq:Martin-Regalado} and study
their stability. While in a physical system injection locking is
possible only at a stable equilibrium, understanding the unstable
equilibrium points provides important insight about the phase space of
the system.

In our previous work~\cite{vonLerber2019alloptical} we concluded that
in the case of linear polarization, a stably injection-locked laser
approximates normalization operation that can be used for arithmetic
computations. In this paper, we widen the scope and explore the
equilibrium points in greater detail. Our main results regarding the
dynamics of system~\eqref{eq:Martin-Regalado} are:
\begin{enumerate}[(i)]
\item If the injected field $u\in\C^2$ is sufficiently weak (small in
  magnitude), then system~\eqref{eq:Martin-Regalado} has nine
  equilibrium points
  (Theorem~\ref{thm:equilibrium-small-dynamics}). If $u$ is
  sufficiently strong (large in magnitude), then it only has a single
  equilibrium point (Theorem~\ref{thm:strong-injection}).
\item Dependence of the equilibrium points on the injected field $u$
  is described in an asymptotic sense in the limits $|u|\to 0$ and
  $|u|\to\infty$ (Theorems~\ref{thm:equilibrium-small-dynamics}
  and~\ref{thm:strong-injection}). A method for calculating the exact
  values of the equilibrium points is provided for weak $u$ in terms
  of an ordinary differential equation
  (Theorem~\ref{thm:solution-from-IVP}).
\item Under the assumption that $\alpha=0$ and that the injected field
  $u$ is weak, it is proved that one of the nine equilibrium points is
  asymptotically stable, while the remaining eight are unstable
  (Theorem~\ref{thm:stability}).
\end{enumerate}

\begin{figure}
  \centering
  \begin{picture}(100,100)
    \put(-50,0){\includegraphics[width=7cm]{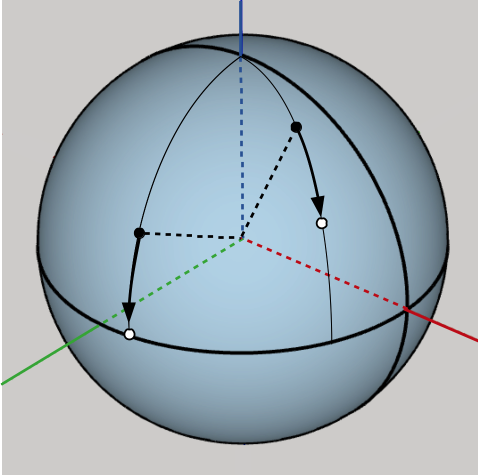}}
    \put(-10,100){a} \put(92,105){b}
  \end{picture}
  \caption{The state of polarization is transformed by the
    injection-locked laser. Schematic illustration on Poincar\'e
    sphere~\cite{shurcliff}: {(a)} A weak injected arbitrary state of
    elliptical polarization ($\bullet$) is projected on equator
    ($\circ$) by the injection-locked laser emission. {(b)} Under a
    strong elliptical state of polarization input, the state of
    polarization of the injection-locked output emission is shifted
    toward the equator, yet, will not reach it.}
  \label{fig:Poincare_sphere}
\end{figure}

The consequence of the aforementioned results is that under weak
injection of elliptically polarized light the injection-locked laser
will emit linearly polarized output such that the input state of
polarization is projected to a linear state of polarization (see
Figure~\ref{fig:Poincare_sphere}a). Under strong injection of
elliptically polarized light, the injection-locked laser will emit
light with an elliptical state of polarization, yet, the polarization
is shifted toward a linear state of polarization, as shown in
Figure~\ref{fig:Poincare_sphere}b.

In the last section of this paper, we will investigate a possibility
to use lasers as nodes of an optical neural network. In general,
optical technologies are commonly used for linear operations, such as
Fourier transformation and matrix multiplications, which come
virtually free by use of lenses, mirrors, and other common light
transforming elements. In this respect, optical solutions have been
proposed for matrix multiplications in optical neural
networks~\cite{shen2017deep, harris2018linear}. However, a neural
network consisting of linear transformations only is impossible, as
such a network is itself linear. As recognized by the optics
community, the nonlinear functions are difficult to realize in
practice, as noted in recent publication
\begin{quote}
  \emph{Despite these positive results, the scheme faces major
    challenges. [...] Then there is the question of the nonlinear
    operation needed to link one set of [Mach-Zehnder Interferometers]
    with another, which [was] simply simulated using a normal
    computer.}~\cite{cartlidge2020optical}
\end{quote}
In this respect, we propose that a laser could provide a useful
nonlinearity. More specifically, a nonlinear activation function of a
node is provided by injection locking; a laser nonlinearly transforms
an injected field (input) into an injection-locked emitted field
(output). As the fields are complex-valued, this also leads in a
natural way to a \emph{complex-valued neural network}.

Complex-valued neural networks are a less studied object than their
real counterpart, nevertheless, they have attracted a considerable
amount of research~\cite{hirose2012complex, aizenberg2011complex,
  hirose2014guest}. A desired quality of any class of neural networks
is the \emph{universal approximation property}, namely, that any
continuous function can be approximated to any degree of accuracy by a
network from that class. For real-valued neural networks, necessary
and sufficient conditions for an activation function to generate a
class of neural networks with the universal approximation property are
known~\cite{leshno1993multilayer,hornik1993some}, and also
quantitative bounds for the approximation exist~\cite{mhaskar,
  yarotsky2017error}. Besides for the theoretical expressiveness of
neural networks, the choice of an activation function affects their
empirical performance, as, among others, it affects the efficacy of
the training
algorithms~\cite{MR3617773}. In~\cite{2017arXiv170907900V} we
considered universality of laser based neural networks with a
complex-valued activation function.

The recent \emph{universal approximation theorem} for complex-valued
neural networks by F.~Voigtlaender~\cite{voigtlaender2020universal}
characterizes those activation functions for which the associated
complex-valued neural networks have the universal approximation
property. In this theorem, the activation function is required to be
defined globally on the complex plane. As the activation function
induced by injection locking is defined only locally in a neighborhood
of the origin, we extend Voigtlaender's theorem by proving a local
version of the universal approximation theorem (Theorem~\ref{thm:UAT}
stated in the Appendix). This theorem and the results about dynamics
of system~\eqref{eq:Martin-Regalado} will prove the following:
\begin{enumerate}
\item[] The class of complex-valued optical neural networks with nodes
  composed of optically injected semiconductor lasers and an
  activation function based on injection locking has the universal
  approximation property, namely, it can approximate any
  complex-valued continuous function to any degree of accuracy
  (Theorem~\ref{thm:ONN}).
\end{enumerate}

The paper is organized as follows. In Sections~\ref{sec:weak-fields}
and~\ref{sec:stability} we assume that the injected field $u$ is weak
and consider equilibrium points of system~\eqref{eq:Martin-Regalado}
and their stability, respectively. In Section~\ref{sec:strong-fields}
we consider the case of a strong injected field. In
Section~\ref{sec:neural-network} we propose a design for an optical
neural network with working principle based on injection locking,
provide a mathematical model for such a network, and prove that these
networks have the universal approximation property. In the Appendix,
we prove a local version of the universal approximation theorem for
complex-valued neural networks.

%
\section{Analysis of equilibrium points and their stability}

\subsection{Equilibrium points with weak injected fields}
\label{sec:weak-fields}

In this section, we study equilibrium points of
system~\eqref{eq:Martin-Regalado} (i.e., points $(E_\pm, N, n)$ at
which the right-hand side of~\eqref{eq:Martin-Regalado} vanishes)
under the assumption that the injected field $u$ is weak and constant
in time. Specifically, we consider injected fields $u$ of the form
\begin{equation}
  \label{eq:lambda-uhat}
  u = \lambda\uhat,  
\end{equation}
where $\uhat\in\C^2\setminus\{0\}$ is fixed and
$\lambda\in\C\setminus\{0\}$ is a small parameter, and we are
interested in the behavior of the equilibrium points as a function of
the parameter $\lambda$.

We assume without loss of generality that $\eta=1$, as this constant
can be incorporated in the injected field $u$. Then we can write
system~\eqref{eq:Martin-Regalado} in an equivalent form
\begin{subequations}
  \label{eq:system}
  \begin{align}
    \frac{d}{dt} E(t) &= -\kappa\big( (1+i\alpha)\, X(N(t),n(t)) E(t) - u \big),\\
    \frac{d}{dt}
    \begin{bmatrix}
      N(t) \\n(t)
    \end{bmatrix}
                      &= -\gamma\left(Y(E(t))
                        \begin{bmatrix}
                          N(t)\\n(t)
                        \end{bmatrix}
    -
    \begin{bmatrix}
      \mu\\0
    \end{bmatrix}\right),
  \end{align}
\end{subequations}
where $E(t) = (E_-(t),E_+(t))$ is a $\C^2$-valued function, and $X$
and $Y$ are matrix-valued functions defined for a vector
$z=(z_1,z_2)\in\C^2$ by
\begin{subequations}
  \label{eq:XY}
  \begin{align}
    X(z) &:=
           \begin{bmatrix}
             1 - (z_1-z_2) & 0\\
             0 & 1 - (z_1+z_2)
           \end{bmatrix},\label{eq:X}\\
    Y(z) &:=
           \begin{bmatrix}
             1+|z|^2 & |z_2|^2 - |z_1|^2\\
             |z_2|^2 - |z_1|^2 & \delta + |z|^2
           \end{bmatrix}\label{eq:Y}
  \end{align}
\end{subequations}
(we use everywhere $(a_1,\ldots,a_n)$ as an alternative notation for a
column vector $\begin{bmatrix}a_1&\cdots&a_n\end{bmatrix}^T$). Above
$|\cdot|$ denotes the absolute value on $\C$ and norm on $\C^2$, and
$\delta := \gamma_s/\gamma>0$ is a dimensionless parameter. The
parameters satisfy ${\delta, \gamma, \kappa}\in(0,\infty)$,
$\alpha\in\R$, and $\mu>1$, and throughout this paper we take them to
be fixed, so that various constants explicit or implicit (as in the
little $o$-notation) in the equations below may depend on them.

\begin{figure}[t!]         
  \centering \input{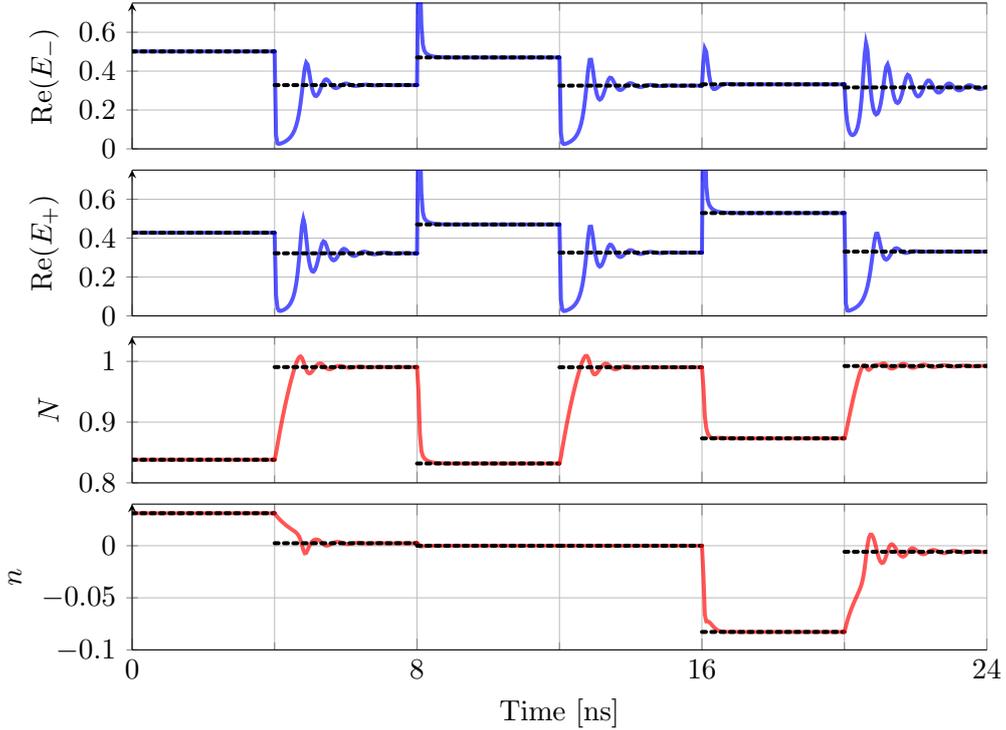}
  \caption[ODE-solution]{Time evolution of the slowly varying
    amplitude $E(t)$ (in circularly polarized basis, blue lines) of an
    electric field emitted by a laser in a case where the slowly
    varying amplitude of an external electric field injected into the
    laser is piecewise constant in time, and corresponding time
    evolution of the parameters $N(t)$ and $n(t)$ (red lines) of the
    laser.
    
    The zero initial value at $t = -4$~ns was used, yet the solution
    is plotted only for $t\ge 0$. In this figure, the injected field
    $u(t) = \lambda(t)\uhat(t)$ has been chosen so that
    $\im(E_\pm(t)) = 0$ for real-valued initial values. Here
    $\lambda(t)=0.25$ for $t\in[-4,0]$ and $t\in[8k, 4(2k+1))$,
    $k\in\{0,1,2\}$, and $\lambda(t)=0.01$ otherwise, and
    $\uhat(t) = \sqrt{\mu-1}\, (\cos\theta(t),\sin\theta(t))$, where
    $\theta(t)=\pi/6$ (corresponding to elliptical polarization) for
    $t\in[-4, 8)$, $\theta(t)=\pi/4$ (linear polarization) for
    $t\in[8, 16)$, and $\theta(t)=11\pi/24$ (nearly circular
    polarization) for $t\in[16, 24)$. After every change in the
    injected field $u$, the laser is seen to quickly stabilize at a
    new equilibrium point. Black dotted lines correspond to the stable
    equilibrium point $E^{(+\textsc{x})}_{\uhat(t)}(\lambda(t))$ (cf.\
    Theorems~\ref{thm:equilibrium-small-dynamics}
    and~\ref{thm:stability}), i.e., they show values of $E$, $N$, and
    $n$ of the laser after a successful injection locking.
    
    In this figure, $\kappa=300$~ns$\inv$, $\mu=1.2$, $\alpha=0$,
    $\gamma=1$~ns$\inv$, and $\delta=\gamma_s/\gamma=1.4$.}
  \label{fig:ODE-solution}
\end{figure}

Figure~\ref{fig:ODE-solution} shows an example of a solution to
system~\eqref{eq:system} with an injected field $u$ that is piecewise
constant.\footnote{All numerical calculations in this article were
  done with Julia~\cite{Julia-2017}. In Figure~\ref{fig:ODE-solution}
  the suite
  DifferentialEquations.jl~\cite{rackauckas2017differentialequations}
  was used.}  After every abrupt change of the injected field $u$, the
solution is seen to quickly settle at a new value (an equilibrium
point of the system).

\begin{proposition}
  \label{prop:system-solution}
  For every initial value $(E_0,N_0,n_0)\in\C^2\times\R\times\R$,
  there exists a unique maximal solution (i.e., a solution that has no
  proper extension that is also a solution) to
  system~\eqref{eq:system} satisfying the initial value at $t=0$. The
  solution is global in forward time, that is, its domain includes
  $[0,\infty)$.
\end{proposition}
\begin{proof}
  A straightforward calculation shows that the right-hand side of
  system~\eqref{eq:system} is locally Lipschitz, which implies that
  for any given initial value, there exists a unique maximal solution
  satisfying the value at $t=0$.
  
  Consider an arbitrary maximal solution
  $(E,N,n):I\to\C^2\times\R\times\R$, where $0\in I\subset\R$, and for
  the sake of a contradiction assume that $[0,\infty)\not\subset
  I$. If $\omega\in\R$ denotes the right endpoint of $I$, then
  $\omega\notin I$ and either
  \begin{equation}\label{eq:blow-up}
    \lim_{
      \substack{
        t \to \omega,\\ t\in I}}
    |E(t)| = \infty\text{ or }
    \lim_{
      \substack{
        t \to \omega,\\ t\in I}}
    \big|\big(N(t),n(t)\big)\big|=\infty
  \end{equation}
  (see \cite[Theorem~{7.6}]{MR1071170}).

  Denote $\nu(t):=(N(t),n(t))\in\R^2$. The function $Y$ is uniformly
  bounded from below, in the sense that there exists $c>0$ such that
  for every $z\in\C^2$ and $y\in\R^2$ it holds that
  \begin{equation}
    \label{eq:Y-lower-bound}
    Y(z)y\cdot y \ge c|y|^2.
  \end{equation}
  With~\eqref{eq:Y-lower-bound} we can estimate
  \begin{equation*}
    \begin{split}
      \frac{1}{2}\Big(\frac{d}{dt}|\nu|^2\Big)(t)
      &=  \dot{\nu}(t)\cdot\nu(t)\\
      &= \gamma\big(-Y(E(t))\nu(t)\cdot\nu(t) + \mu N(t)\big)\\
      &\le C_1(1+|\nu(t)|^2),
    \end{split}
  \end{equation*}
  where $0<t<\omega$ and $C_1\in\R$ is a constant. This inequality
  together with Gr\"onwall's lemma yields $|\nu(t)| \le C_2$ for every
  $0\le t<\omega$, where $C_2\ge 0$ is another constant.

  The fact that $\nu$ is bounded on $[0,\omega)$ implies that the
  function $t\mapsto X(\nu(t))$ is also bounded there. Then similar
  reasoning as above (involving Gr\"onwall's lemma) shows that $E$ is
  bounded on $[0,\omega)$.  This contradicts with~\eqref{eq:blow-up},
  and therefore $[0,\infty)\subset I$.
\end{proof}

Following theorem is the main result of this section. Its essential
content is that with sufficiently weak injected fields of the form
$u=\lambda\uhat$ system~\eqref{eq:system} has nine distinct
equilibrium points, and that the equilibrium points depend
continuously on $\lambda\in\C\setminus\{0\}$ with asymptotics given
by~\eqref{eq:E-approximations}. In the statement of the theorem, the
requirement that $\uhat_-\neq 0$ and $\uhat_+\neq 0$ means physically
that the field is not \emph{circularly polarized}, while
$|\uhat_-|=|\uhat_+|$ means that the field is \emph{linearly
  polarized}. The function $y:\R^2\to\R^2$ is defined by
\begin{align}
  y(x)  &:= Y(x)\inv
          \begin{bmatrix}
            \mu\\0
          \end{bmatrix}
  = \frac{\mu}{\det Y(x)}
  \begin{bmatrix}
    \delta+|x|^2\\
    x_1^2-x_2^2
  \end{bmatrix}, \text{ where} \label{eq:def-y}\\
  \det Y(x) &= \delta + (1+\delta)|x|^2+4x_1^2x_2^2>0\label{eq:detY}
\end{align}
(the function $Y$ is defined in~\eqref{eq:Y}).

\begin{theorem}
  \label{thm:equilibrium-small-dynamics}
  Consider injected external field with amplitude $\lambda\uhat$,
  where $\lambda\in\C$ and $\uhat=(\uhat_-, \uhat_+)\in\C^2$ satisfies
  $\uhat_-\neq 0$ and $\uhat_+\neq 0$. There exists a constant
  $\ell=\ell(\uhat)>0$ and a family $\{E_{\uhat}^{(j)}\}_{j\in\J}$ of
  nine continuous functions
  \begin{equation}
    \label{eq:E-lambda-j}
    E_{\uhat}^{(j)} : \{\lambda\in\C : 0<|\lambda|<\ell\}\to\C^2,\,
    j\in\J
    :=\{\textsc{0}, \pm\textsc{l},\pm\textsc{r}, \pm\textsc{x}, \pm\textsc{y}\},
  \end{equation}
  with pairwise distinct values that have the following properties:
  \begin{enumerate}[(i)]
  \item If in system~\eqref{eq:system} the injected field is of the
    form $u=\lambda\uhat$ with $0<|\lambda|<\ell$, then a triple
    $(E, N, n)\in\C^2\times\R\times\R$ is an equilibrium point (a
    time-independent solution) of the system, if and only if
    \begin{equation*}
      E = E_{\uhat}^{(j)}(\lambda)\text{ for some } j\in\J,
      \text{ and } (N,n) = y(|E_-|,|E_+|). 
    \end{equation*}
  \item The functions $E_{\uhat}^{(j)}$ have following asymptotics as
    $\lambda\to 0$:
    \begin{subequations}
      \label{eq:E-approximations}
      \begin{align}
        E_{\uhat}^{(\textsc{0})}(\lambda)
        &= e^{i\theta}\frac{\lambda}{|\lambda|}
          \left( |\lambda| \what^{(\textsc{0})}
          + o(\lambda) \right),\label{eq:E0-approximation}\\
        E_{\uhat}^{(\pm\textsc{l})}(\lambda)
        &= e^{i\theta}\frac{\lambda}{|\lambda|}
          \left( \pm\sqrt{\frac{\delta(\mu-1)}{1+\delta}}
          \begin{bmatrix}
            \uhat_-/|\uhat_-| \\
            0
          \end{bmatrix}
        + |\lambda|\,\,\what^{(\textsc{l})}
        + o(\lambda) \right),\\
        E_{\uhat}^{(\pm\textsc{r})}(\lambda)
        &= e^{i\theta}\frac{\lambda}{|\lambda|}
          \left( \pm\sqrt{\frac{\delta(\mu-1)}{1+\delta}}
          \begin{bmatrix}
            0 \\
            \uhat_+/|\uhat_+|
          \end{bmatrix}
        + |\lambda|\,\what^{(\textsc{r})}
        + o(\lambda)\right),\\
        E_{\uhat}^{(\pm\textsc{x})}(\lambda)
        &= e^{i\theta}\frac{\lambda}{|\lambda|}
          \left( \pm\sqrt{\frac{\mu-1}{2}}
          \begin{bmatrix}
            \uhat_-/|\uhat_-|\\
            \uhat_+/|\uhat_+|
          \end{bmatrix}
        +|\lambda|\,\what^{(\textsc{x})}
        + o(\lambda)\right),\\
        E_{\uhat}^{(\pm\textsc{y})}(\lambda)
        &= e^{i\theta}\frac{\lambda}{|\lambda|}
          \left(\pm\sqrt{\frac{\mu-1}{2}}
          \begin{bmatrix}
            \phantom{-}\uhat_-/|\uhat_-|\\
            -\uhat_+/|\uhat_+|
          \end{bmatrix}
        +|\lambda|\,\what^{(\textsc{y})}
        + o(\lambda)\right),
      \end{align}
    \end{subequations}
    where $\theta := -\arg(1+i\alpha)$ and
    \begin{align*}
      \what^{(\textsc{0})}
      & := \frac{-1}{|1+i\alpha|(\mu-1)} \uhat\\
      \what^{(\textsc{l})}
      &:= \frac{1}{2|1+i\alpha|(\mu-1)}
        \begin{bmatrix*}[r]
          \mu\,\uhat_-\\
          -(1+\delta)\,\uhat_+
        \end{bmatrix*},\\
      \what^{(\textsc{r})}
      &:= \frac{1}{2|1+i\alpha|(\mu-1)}
        \begin{bmatrix*}[r]
          -(1+\delta)\,\uhat_-\\
          \mu\,\uhat_+
        \end{bmatrix*},\\
      \what^{(\textsc{x})}
      &:= \frac{1}{4|1+i\alpha|(\mu-1)}
        \begin{bmatrix*}[r]
          \big(2\mu+\delta-1 + (1-\delta)|\uhat_+|/|\uhat_-|\big)\,\uhat_- \\
          \big((1-\delta)|\uhat_-|/|\uhat_+| +
          2\mu+\delta-1\big)\,\uhat_+
        \end{bmatrix*},\\
      \what^{(\textsc{y})}
      &:= \frac{1}{4|1+i\alpha|(\mu-1)}
        \begin{bmatrix*}[r]
          \big(2\mu+\delta-1 + (\delta-1)|\uhat_+|/|\uhat_-|\big)\,\uhat_- \\
          \big((\delta-1)|\uhat_-|/|\uhat_+| +
          2\mu+\delta-1\big)\,\uhat_+
        \end{bmatrix*}.\\
    \end{align*}
  \item Furthermore, if $|\uhat_-|=|\uhat_+|$ and
    $j\in\{0,\pm\textsc{x}\}$, then for every $\lambda$ with
    $0<|\lambda|<\ell$ it holds that
    \begin{equation*}
      E_{\uhat}^{(j)}(\lambda) = \rho^{(j)}(\lambda)\uhat
    \end{equation*}
    for some $\rho^{(j)}(\lambda)\in\C$.
  \end{enumerate}
\end{theorem}

\begin{remark}
  As $\lambda\to 0$, the amplitude $E_{\uhat}^{(\textsc{0})}(\lambda)$
  vanishes, the amplitudes $E_{\uhat}^{(\pm\textsc{l})}(\lambda)$ and
  $E_{\uhat}^{(\pm\textsc{r})}(\lambda)$ become left and right
  circularly polarized, respectively, and the amplitudes
  $E_{\uhat}^{(\pm\textsc{x})}(\lambda)$ and
  $E_{\uhat}^{(\pm\textsc{y})}(\lambda)$ become linearly polarized and
  orthogonal to each other. The index set $\J$ is chosen to reflect
  this fact. Note that as $\lambda\to 0$, on the normalized Poincar\'e
  sphere the amplitudes $E_{\uhat}^{(\pm\textsc{x})}(\lambda)$
  approach the projection of $\uhat$ onto the equator, and the
  amplitudes $E_{\uhat}^{(\pm\textsc{y})}(\lambda)$ approach the
  antipodal point of that projection.
\end{remark}

\begin{remark}
  At the expense of a more complicated statement, the theorem can be
  modified to hold also in the case $\uhat_-=0$ or $\uhat_+=0$. The
  reason why this case is special is that if a point $(E_-,E_+,N,n)$
  is an equilibrium point of system~\eqref{eq:system} with injected
  field (say) $u=(0,\lambda\uhat_+)$, then for every $\phi\in\R$ the
  point $(e^ {i\phi}\,E_-,E_+, N, n)$ is an equilibrium point of the
  system, also. Thus, instead of distinct equilibrium points, there
  will be disjoint sets of equilibrium points. See also
  Proposition~\ref{prop:algebraic-version},
  Remark~\ref{remark:non-uniqueness}, and
  Theorem~\ref{thm:solution-from-IVP} below.
\end{remark}

\begin{figure}[tp]
  \centering \input{ax_E.tikz}
  \caption[Equilibrium dynamics]{Values on the $\re(E_\pm)$-plane
    (circularly polarized basis) of the slowly varying amplitude $E$
    of an electric field of a laser at the stable and unstable
    equilibrium points as the magnitude $\lambda$ of an external
    optical injection $u=\lambda\uhat$ varies.
    
    For small $|\lambda|$ the laser has nine equilibrium points
    (Theorem~\ref{thm:equilibrium-small-dynamics}). Solid lines denote
    paths traced by real parts of the points when
    $\uhat = \sqrt{\mu-1}(\cos\theta, \sin\theta)$ and
    $\theta = \pi/4$ (linear polarization), and
    $\lambda\in[-1/4, 1/4]$ varies. In this figure, $\uhat$ has been
    chosen so that the equilibrium points are real-valued for
    $\lambda\in\R$ and so that the intensity of the injected field
    $u=\lambda\uhat$ at $\lambda=1$ is equal to the intensity of the
    emitted field $E$ of the free-running laser.
    
    As $\lambda$ increases, the points move in the directions indicated
    by the arrows.  At $\lambda = -1/4$ only one of the points exists
    (it is located at {(i)}). As $\lambda$ increases, eight new points
    appear. First at $\lambda\approx -0.072$ two points appear at {(a)}
    and start moving in opposite directions. At $\lambda\approx -0.071$
    one of these points has moved to {(b)}, where it splits into three.
    At $\lambda\approx -0.057$ two points appear at each {(c)}. The
    circled dots denote locations of the points at $\lambda=0$. As
    $\lambda$ grows, eight of the points disappear (at {(d)}
    ($\lambda\approx 0.057$), {(e)} ($\lambda\approx 0.71$), and {(f)}
    ($\lambda\approx 0.072$)). The paths were calculated from the
    functions $h^{(j)}_{\rhat}$ (cf. Figure~\ref{fig:h-solution})
    via~\eqref{eq:E-N-n-u} and~\eqref{eq:u-from-r-phi}.
    
    For $-1/4\le\lambda<0$ only the equilibrium point on the path from
    {(i)} to {(ii)} is stable, for $0<\lambda\le 1/4$ the same is
    true for the equilibrium point on the path from {(iii)} to {(iv)}
    (cf.\ Figure~\ref{fig:stability}). Consequently, at $\lambda=0$,
    the unique stable equilibrium point of the system jumps from
    {(ii)} to {(iii)}.
    
    The parameters $\kappa$, $\mu$, $\alpha$, $\gamma$, and $\delta$
    are those of Figure~\ref{fig:ODE-solution}. The dotted and dashed
    paths are interpreted analogously. In these paths
    $\theta\in\{\pi/6, 11\pi/24\}$ (elliptical polarizations).}
  \label{fig:paths}
\end{figure}

Figure~\ref{fig:paths} shows values of the nine equilibrium points
from Theorem~\ref{thm:equilibrium-small-dynamics} as the magnitude
$\lambda\in\R$ of an external optical injection $u=\lambda\uhat$
varies. In the dimensionless units of system~\eqref{eq:system} the
intensity of the free running laser, i.e., $|E|^2$ at a stable
equilibrium point of~\eqref{eq:system} when $u=0$, is $\mu-1$. In the
figure $\uhat$ has been chosen so that at $|\lambda|=1$ the intensity
$|u|^2=|\uhat|^2$ of the external injected field is also $\mu-1$. For
the laser parameters used in the figure, the injected field is
sufficiently weak in the sense of
Theorem~\ref{thm:equilibrium-small-dynamics}, namely, in the sense
that the nine equilibrium points of the theorem exist, if
$|\lambda| < 0.057$, i.e., if the injected field does not exceed in
magnitude 5.7~\% of the emitted field of the free running laser. In
practice this value would depend also on experimental setup details
such as the coupling efficiency.

As a real-valued amplitude $E=(E_-,E_+)\in\R^2$ is linearly polarized
if and only if $E_-=\pm E_+$, it is seen from Figure~\ref{fig:paths}
that even if the injected field is linearly polarized, only three of
the nine equilibrium points have a linear state of polarization, while
the remaining six equilibrium points have an elliptical state of
polarization.

We prove Theorem~\ref{thm:equilibrium-small-dynamics} at the end of
this section after developing some preliminary results. We begin by
transforming the problem of finding equilibrium points of
system~\eqref{eq:system} from $\C^2\times\R\times\R$ into a problem of
finding solutions from $\R^2$ to a system of two bivariate
polynomials:
\begin{proposition}
  \label{prop:algebraic-version}
  Let $X$ and $y$ be the functions defined in~\eqref{eq:X}
  and~\eqref{eq:def-y}.
  \begin{enumerate}[(i)]
  \item Fix a vector $r=(r_1,r_2)\in[0,\infty)\times[0,\infty)$, and
    suppose $x=(x_1,x_2)\in\R^2$ satisfies
    \begin{equation}
      \label{eq:X-y-x-x-r}
      X(y(x))x = r.
    \end{equation}
    Let $\phi_\pm\in\R$, and define a vector $E\in\C^2$ and numbers
    ${N,n}\in\R$ by
    \begin{equation}
      \label{eq:E-N-n-u}
      E :=
      \begin{bmatrix}
        x_1\,e^{i\phi_-}\\
        x_2\,e^{i\phi_+}
      \end{bmatrix}\text{ and }
      \begin{bmatrix}
        N\\n
      \end{bmatrix}
      := y(x).
    \end{equation}
    Then the triple $(E,N,n)$ is an an equilibrium point of
    system~\eqref{eq:system} with the injected electric field
    \begin{equation}
      \label{eq:u-from-r-phi}
      u := (1+i\alpha)
      \begin{bmatrix}
        r_1\,e^{i\phi_-}\\
        r_2\,e^{i\phi_+}
      \end{bmatrix}.
    \end{equation}
  \item Suppose a triple $(E,N,n)\in\C^2\times\R\times\R$ is an
    equilibrium point of system~\eqref{eq:system} with some injected
    electric field $u\in\C^2$. Then there exists numbers
    $\phi_\pm\in\R$ and vectors $r\in[0,\infty)\times[0,\infty)$ and
    $x\in\R^2$ such that equations~\eqref{eq:X-y-x-x-r}
    to~\eqref{eq:u-from-r-phi} hold.
  \end{enumerate}
\end{proposition}

\begin{remark}
  \label{remark:non-uniqueness}
  An arbitrary field $u=(u_-,u_+)\in\C^2$ uniquely determines the
  numbers $r_j\ge 0$ in~\eqref{eq:u-from-r-phi}. If $u_-\neq 0$ and
  $u_+\neq 0$, then also the numbers $e^{i\phi_\pm}$ are uniquely
  determined, and therefore a solution $x\in\R^2$
  of~\eqref{eq:X-y-x-x-r} corresponds via~\eqref{eq:E-N-n-u} to a
  unique equilibrium point of system~\eqref{eq:system}. But if (say)
  $u_-=0$ and $x$ is a solution of~\eqref{eq:X-y-x-x-r} with
  $x_1\neq 0$, then there exists a continuum of equilibrium points of
  system~\eqref{eq:system} corresponding to $x$ due to the arbitrary
  choice of $\phi_-\in\R$ in~\eqref{eq:u-from-r-phi}.
\end{remark}

\begin{proof}[Proof of Proposition~\ref{prop:algebraic-version}]
  For a vector $\phi=(\phi_-,\phi_+)\in\R^2$ denote
  \begin{equation*}
    J_\phi :=
    \begin{bmatrix}
      e^{i\phi_-}&0\\
      0&e^{i\phi_+}
    \end{bmatrix}
    \in\C^{2\times 2}.
  \end{equation*}
  Then for every $z\in\C^2$ the matrices $X(z)$ and $J_\phi$ commute,
  and $Y(J_\phi z)=Y(z)$.

  For proving the first part of the proposition assume that $x\in\R^2$
  and $r\in[0,\infty)\times[0,\infty)$ satisfy~\eqref{eq:X-y-x-x-r},
  and let $E = J_\phi x$, $(N,n)=y(x)$, and $u=(1+i\alpha)J_\phi r$ be
  as in~\eqref{eq:E-N-n-u} and~\eqref{eq:u-from-r-phi}. Then
  \begin{align*}
    -\kappa\big((1+i\alpha)X(N,n)E-u\big)
    &= -\kappa(1+i\alpha)J_\phi\big(X(y(x))x-r\big)
      = 0,\text{ and}\\
    -\gamma\left(Y(E)
    \begin{bmatrix}
      N\\n
    \end{bmatrix}
    -
    \begin{bmatrix}
      \mu\\0
    \end{bmatrix}\right)
    &= -\gamma\left(Y(x)Y(x)^{-1}
      \begin{bmatrix}
        \mu\\0
      \end{bmatrix}
    -
    \begin{bmatrix}
      \mu\\0
    \end{bmatrix}
    \right)
    =0,
  \end{align*}
  so the point $(E,N,n)$ is an equilibrium point of
  system~\eqref{eq:system} with the injected electric field $u$.

  For proving the second part of the proposition assume a point
  $(E,N,n)\in\C^2\times\R\times\R$ is an equilibrium point of
  system~\eqref{eq:system} with injected electric field $u\in\C^2$,
  and find vectors $x\in\R^2$ and $\phi=(\phi_-,\phi_+)\in\R^2$ such
  that $E=J_\phi x$ and
  \begin{equation}
    \label{eq:re-is-pos}
    \re\left(\frac{e^{-i\phi_\pm}u_\pm}{1+i\alpha}\right) \ge 0.
  \end{equation}
  Then from above and the definition of an equilibrium point it
  follows that
  \begin{equation*}
    \begin{bmatrix}
      \mu\\0
    \end{bmatrix}
    = Y(E)
    \begin{bmatrix}
      N\\n
    \end{bmatrix}
    = Y(x)
    \begin{bmatrix}
      N\\n
    \end{bmatrix}
    \text{ and }
    u = (1+i\alpha)X(N,n)J_\phi x.
  \end{equation*}
  This implies $(N,n)=y(x)$, and consequently
  $u = (1+i\alpha)J_\phi X(y(x))x$.

  Now define $r:=X(y(x))x\in\R^2$. Then it only remains to show that
  $r_j\ge 0$, but this follows from~\eqref{eq:re-is-pos}, since
  $r=(1+i\alpha)^{-1}J_{-\phi}u$.
\end{proof}
   
\begin{proposition}
  \label{prop:zeros-at-0}
  A vector $x\in\R^2$ satisfies $X(y(x))x=0$, if and only if
  $x = x^{(j)}$ for some $j\in\J$ (the index set $\J$ is defined
  in~\eqref{eq:E-lambda-j}), where
  \begin{subequations}
    \label{eq:xk}
    \begin{align}
      x^{(\textsc{0})} &:=
                         \begin{bmatrix}
                           0\\0
                         \end{bmatrix},\\
      x^{(\pm\textsc{l})} &:= \pm\sqrt{\frac{\delta(\mu-1)}{1+\delta}}
                            \begin{bmatrix}
                              1\\0
                            \end{bmatrix},\\
      x^{(\pm\textsc{r})} &:= \pm\sqrt{\frac{\delta(\mu-1)}{1+\delta}}
                            \begin{bmatrix}
                              0\\1
                            \end{bmatrix},\\
      x^{(\pm\textsc{x})} &:= \pm\sqrt{\frac{\mu-1}{2}}
                            \begin{bmatrix}
                              1\\1
                            \end{bmatrix},\\
      x^{(\pm\textsc{y})} &:= \pm\sqrt{\frac{\mu-1}{2}}
                            \begin{bmatrix}
                              \phantom{-}1\\-1
                            \end{bmatrix}.
    \end{align}
  \end{subequations}
\end{proposition}
\begin{proof}
  Suppose that $X(y(x))x=0$, or equivalently that
  \begin{equation}
    \label{eq:y-x-equality}
    y_2(x)
    \begin{bmatrix}
      1 & \phantom{-}0 \\ 0 & -1
    \end{bmatrix}
    x = (y_1(x)-1)x,
  \end{equation}
  where $y(x) = (y_1(x), y_2(x))$. It follows that if
  $x\neq x^{(\textsc{0})}$, then $|y_2(x)|=|y_1(x)-1|$.

  Consider first the case $y_1(x)-1=y_2(x)\neq
  0$. Then~\eqref{eq:y-x-equality} implies that $x$ is of the form
  $(c,0)$ for some $c\in\R$. To find the possible values of $c$,
  insert the candidate vector into $y_1(x)-1=y_2(x)$ and solve for
  $c$. This shows that
  $x\in \{{x^{(+\textsc{l})},x^{(-\textsc{l})}}\}$.

  If $1-y_1(x)=y_2(x)\neq 0$, an analogous reasoning shows that then
  $x\in \{{x^{(+\textsc{r})},x^{(-\textsc{r})}}\}$.

  Consider the last case, namely $y_2(x)=0$ and $y_1(x)=1$. Then
  $x_1^2=x_2^2$, and inserting $x=(x_1,\pm x_1)$ into $y_1(x)=1$ and
  solving for $x_1$ shows that $x_1^2=x_2^2 = (\mu-1)/2$. Taking into
  account all possible sign combinations yields
  $x\in\{
  {x^{(+\textsc{x})},x^{(-\textsc{x})},x^{(+\textsc{y})},x^{(-\textsc{y})}}
  \}$.

  On the other hand, a direct calculation shows that
  $X(y(x^{(j)}))x^{(j)}=0$ for every $j\in\J$.
\end{proof}

Fix nonzero $\rhat=(\rhat_1,\rhat_2)\in\R^2$ and define a function
\begin{equation}
  \label{eq:def-F}
  F_{\rhat}(s,x) := X(y(x))x-s\rhat \qquad(s\in\R, \, x\in\R^2). 
\end{equation}
Our plan is to first find all zeros of $F_{\rhat}(s,\cdot)$ for small
$s$, and then, assuming that the injected field $u$ in
system~\eqref{eq:system} is sufficiently weak, with
Proposition~\ref{prop:algebraic-version} convert these zeros to
equilibrium points of the system.

The Jacobian matrix of $F_{\rhat}$ with respect to $x$ will be denoted
by $D_x F_{\rhat}(x)$ (as the Jacobian is independent of $s$, it is
suppressed from the notation). A straightforward calculation shows
that
\begin{equation}
  \label{eq:DxF-expression}
  D_x F_{\rhat}(x)
  = I_2 + \frac{1}{\det Y(x)^2}
  \begin{bmatrix}
    p_{11}(x) & p_{12}(x)\\
    p_{21}(x) & p_{22}(x)
  \end{bmatrix},
\end{equation}
where $I_2\in\R^{2\times 2}$ is the identity matrix,
\begin{align*}
  p_{11}(x_1,x_2) &:= \mu(\delta+2x_2^2)(-\delta+(1+\delta)(x_1^2-x_2^2)+4x_1^2x_2^2),\\
  p_{12}(x_1,x_2) &:= 2\mu(\delta-1)(\delta+2x_1^2)x_1x_2,\\
  p_{21}(x_1,x_2) &:= p_{12}(x_2, x_1),\text{ and}\\
  p_{22}(x_1,x_2) &:= p_{11}(x_2, x_1)
\end{align*}
(an expression for $\det Y(x)$ is given in \eqref{eq:detY}).
 
\begin{proposition}
  \label{prop:F-properties}
  \begin{enumerate}[(i)]
  \item The matrices $D_xF_{\rhat}(x^{(j)})$, $j\in\J$, are
    invertible, and
    \begin{align*}
      \big[D_xF_{\rhat}(x^{(\textsc{0})})\big]\inv
      &= -\frac{1}{\mu-1} I_2,\\
      \big[D_xF_{\rhat}(x^{(\pm\textsc{l})})\big]\inv
      &= \frac{1}{2}\frac{1}{\mu-1}
        \begin{bmatrix}
          \mu & 0\\
          0 & -(1+\delta)
        \end{bmatrix},\\
      \big[D_xF_{\rhat}(x^{(\pm\textsc{r})})\big]\inv%
      &=\frac{1}{2}\frac{1}{\mu-1}
        \begin{bmatrix}
          -(1+\delta) & 0\\
          0 & \mu
        \end{bmatrix},\\
      \big[D_xF_{\rhat}(x^{(\pm\textsc{x})})\big]\inv%
      &=\frac{1}{4}\frac{1}{\mu-1}
        \begin{bmatrix}
          2\mu+\delta-1 & 1-\delta\\
          1-\delta & 2\mu+\delta-1
        \end{bmatrix},\\
      \big[D_xF_{\rhat}(x^{(\pm\textsc{y})})\big]\inv%
      &=\frac{1}{4}\frac{1}{\mu-1}
        \begin{bmatrix}
          2\mu+\delta-1 & \delta-1\\
          \delta-1 & 2\mu+\delta-1
        \end{bmatrix}.
    \end{align*}
  \item For nonzero $x\in\R^2$ denote
    \begin{equation*}
      \xhat := |x|\inv x\text{ and }
      \xhat_\perp := |x|\inv
      \begin{bmatrix}
        \phantom{-}x_2\\
        -x_1
      \end{bmatrix}.
    \end{equation*}
    Then
    \begin{equation*}
      X(y(x))x = |x|(a(x)\xhat + b(x)\xhat_\perp)
      \qquad(x\in\R^2\setminus\{0\}),
    \end{equation*}
    where following estimates hold for the functions
    ${a,b}:\R^2\setminus\{0\}\to\R$:
    \begin{subequations} \label{eq:ab-estimates}
      \begin{align}
        0 \le 1 - a(x) &< \mu\,\min\left\{1, \frac{1}{|x|^2}\right\}
                         \text{, and }\label{eq:a-estimate}\\
        |b(x)| &< \mu\,\min\left\{\frac{1}{1+\delta}, \frac{1}{(1+\delta)^{2/3}|x|^{2/3}}\right\}\label{eq:b-estimate}
      \end{align}
    \end{subequations}
    (recall that $\mu>1$). In particular, $a(x)\to 1$ and $b(x)\to 0$
    as $|x|\to\infty$.
  \end{enumerate}
\end{proposition}
\begin{proof}
  Inserting the value of $x^{(j)}$ from~\eqref{eq:xk} into the
  expression~\eqref{eq:DxF-expression} of $D_xF_{\rhat}$ and inverting
  yields {(i)}.
  
  For {(ii)}, consider a vector $x\in\R^2\setminus\{0\}$. A
  calculation shows that
  \begin{equation*}
    1-a(x)
    =\frac{\mu}{|x|^2}\,\frac{\delta|x|^2+4x_1^2x_2^2}{\delta + (1+\delta)|x|^2+4x_1^2x_2^2}
    \in\left(0, \frac{\mu}{|x|^2}\right).
  \end{equation*}
  On the other hand, above together with the inequality
  $4x_1^2x_2^2/|x|^2\le|x|^2$ yields
  \begin{equation*}
    1-a(x)
    = \mu \frac{\delta+4x_1^2x_2^2/|x|^2}{\delta+(1+\delta)|x|^2+4x_1^2x_2^2}
    < \mu.
  \end{equation*}
  Inequality~\eqref{eq:a-estimate} is now proved.

  Regarding the second inequality, note that
  \begin{equation}
    \label{eq:b-expression}
    |b(x)|
    = 2\mu\frac{|x_1x_2|}{|x|^2}\frac{|x_ 1^2-x_2^2|}{\delta + (1+\delta)|x|^2+4x_1^2x_2^2}.
  \end{equation}
  Let $c\ge 0$ be a parameter and consider two cases: If
  $|x_1x_2| < c|x|/2$, then
  \begin{equation}
    \label{eq:b-1}
    |b(x)| < \mu\frac{c}{|x|}\frac{1}{1+\delta}.
  \end{equation}
  If $|x_1x_2|\ge c|x|/2$, applying the inequality $2|x_1x_2|\le|x|^2$
  to~\eqref{eq:b-expression} shows that
  \begin{equation}
    \label{eq:b-2}
    |b(x)|
    \le\mu\frac{|x_1^2-x_2^2|}{\delta+(1+\delta)|x|^2+c^2|x|^2}
    < \mu\frac{1}{1+\delta+c^2}.
  \end{equation}
  Inequalities~\eqref{eq:b-1} and~\eqref{eq:b-2} hold for every
  $c\ge0$. Choosing $c=0$ yields one part of~\eqref{eq:b-estimate},
  choosing $c=(1+\delta)^{1/3}|x|^{1/3}$ yields the other part.
\end{proof}

\begin{proposition}
  \label{prop:implicit-function-theorem}
  There exists $\ell>0$ and smooth functions
  $h^{(j)}_{\rhat}:(-\ell,\ell)\to\R^2$, $j\in\J$, such that the
  following holds: $h^{(j)}_{\rhat}(\textsc{0})=x^{(j)}$ for every
  $j\in\J$, and if $s\in(-\ell,\ell)$, then
  \begin{equation}
    \label{eq:F-h-iff}
    F_{\rhat}(s,x) = 0,\text{ if and only if } x = h^{(j)}_{\rhat}(s)\text{ for some } j\in\J.
  \end{equation}
  ($F_{\rhat}$ is defined in~\eqref{eq:def-F}, $x^{(j)}$
  in~\eqref{eq:xk}, and $\J$ in~\eqref{eq:E-lambda-j}.)  Furthermore,
  if $\rhat_1=\rhat_2$ and $j\in\{0, \pm\textsc{x}\}$, then
  $h_{\rhat}^{(j)}$ is of the form
  \begin{equation}
    \label{eq:r1r2}
    h^{(j)}_{\rhat}(s) = (\eta^{(j)}(s), \eta^{(j)}(s)) 
  \end{equation}
  for some function $\eta^{(j)}:(-\ell,\ell)\to\R$.
\end{proposition}
\begin{proof}
  Recall that $\rhat\neq 0$ by assumption. By {(i)} of
  Proposition~\ref{prop:F-properties} and the implicit function
  theorem there exists neighborhoods $V^{(j)}\subset\R$ of $0\in\R$
  and $W^{(j)}\subset\R^2$ of $x^{(j)}$ and smooth functions
  $h^{(j)}_{\rhat}:V^{(j)}\to W^{(j)}$ with
  $h^{(j)}_{\rhat}(\textsc{0})=x^{(j)}$ such that $F_{\rhat}(s,x)=0$
  for $(s,x)\in V^{(j)}\times W^{(j)}$, if and only if
  $x=h^{(j)}_{\rhat}(s)$.

  Regarding the other direction of~\eqref{eq:F-h-iff}, it is enough to
  show that there exists $\ell>0$ such that
  \begin{equation}
    \label{eq:subset-condition}
    (-\ell,\ell)\subset\bigcap_{j\in\J} V^{(j)}
  \end{equation}
  and that $F_{\rhat}(s,x)= 0$ implies that either
  $(s,x)\in\bigcup_{j\in\J} V^{(j)}\times W^{(j)}$ or $|s|\ge\ell$.

  If a pair $(s,x)\in\R\times\R^2$ satisfies $F_{\rhat}(s,x)=0$ and
  $|x|>\sqrt{2\mu}$, then by the Pythagorean theorem (with the
  notation of Proposition~\ref{prop:F-properties}) we have
  \begin{equation*}
    s^2|\rhat|^2
    = |X(y(x))x|^2
    = |x|^2(a(x)^2+b(x)^2)
    > \frac{|x|^2}{4},
  \end{equation*}
  where the last inequality holds because $a(x)>1/2$
  by~\eqref{eq:a-estimate}. This implies that
  $|s|>\sqrt{\mu}/(\sqrt{2}|\rhat|)$, which together with the
  continuity of $F_{\rhat}$ shows that the set
  \begin{equation*}
    K:=
    F_{\rhat}\inv(\{0\})
    \cap\left\{(s,x) : |s|\le\sqrt{\mu}/(2|\rhat|)\right\}
    \cap\Big(\bigcup_{j\in\J} V^{(j)}\times W^{(j)}\Big)^\complement
    \subset\R\times\R^2
  \end{equation*}
  is compact.
  
  By Proposition~\ref{prop:zeros-at-0} the set $K$ and the closed set
  $\{0\}\times\R^2$ are disjoint. Let $d>0$ be the distance between
  those sets ($d=+\infty$ if $K=\emptyset$), and consider a pair
  $(s,x)$ such that $|s|\le\sqrt{\mu}/(\sqrt{2}|\rhat|)$ and
  $F_{\rhat}(s,x)=0$. Now if $(s,x)\in K$, then $|s|\ge d$, and if
  $(s,x)\notin K$, then
  $(s,x)\in\bigcup_{j\in\J} V^{(j)}\times W^{(j)}$. Consequently, if
  we choose $\ell>0$ small enough so that~\eqref{eq:subset-condition}
  and $\ell<\min\{d,\sqrt{\mu}/(\sqrt{2}|\rhat|)\}$ hold,
  then~\eqref{eq:F-h-iff} holds for every $|s|<\ell$.

  Finally, if $\rhat_1=\rhat_2$ and $\eta\in\R$, then
  $F_{\rhat}(s, (\eta,\eta)) = 0$, if and only if
  \begin{equation}
    \label{eq:scalar-implicit-fn}
    \eta\left(1-\frac{\mu(\delta+2\eta^2)}{\delta+2(1+\delta)\eta^2+4\eta^4}\right)
    - s\rhat_1 = 0.
  \end{equation}
  The implicit function theorem shows that in some neighborhoods of
  $(0,0)\in\R\times\R$ and $(0,\pm\sqrt{(\mu-1)/2})\in\R\times\R$
  equality~\eqref{eq:scalar-implicit-fn} implicitly defines
  $\eta=\eta(s)$, and consequently, if $j\in\{0, \pm\textsc{x}\}$ and
  $s$ is small enough, then $h^{(j)}(s)=(\eta(s),\eta(s))$.
\end{proof}

Following theorem shows that system~\eqref{eq:system} has at least
nine disjoint families of equilibrium points provided that the
injected field $u$ is weak enough. These families correspond to nine
distinct solutions of $F_{\rhat}(s,\cdot)=0$, where $s>0$ is a fixed
parameter related to the strength of the field $u$. These solutions
can be found by solving an initial value problem for an ordinary
differential equation in $s$. As the initial value problem is easy to
solve numerically, the theorem provides a computational method for
obtaining numerical values for the nine families of equilibrium
points.

\begin{theorem}
  \label{thm:solution-from-IVP}
  Fix $\uhat=(\uhat_-,\uhat_+)\in\C^2$ (with the possibility
  $\uhat_-=0$ or $\uhat_+=0$ allowed), and consider
  system~\eqref{eq:system} with $u=\lambda\uhat$.  Define
  \begin{equation*}
    \rhat := \frac{1}{|1+i\alpha|}
    \begin{bmatrix}
      |\uhat_-|\\
      |\uhat_+|
    \end{bmatrix}
    \in[0,\infty)\times[0,\infty)
  \end{equation*}
  and choose numbers $\phi_\pm\in\R$ such that
  \begin{equation}
    \label{eq:phi-choice}
    \uhat_\pm = |\uhat_\pm|e^{i\phi_\pm}.
  \end{equation}
  Let $y:\R^2\to\R^2$ and $\J$ be as defined in~\eqref{eq:def-y}
  and~\eqref{eq:E-lambda-j}, respectively, and define
  $\theta:=-\arg(1+i\alpha)$.

  Fix $j\in\J$. Suppose $I\subset\R$ is an interval containing the
  origin and
  \begin{equation*}
    h = (h_1, h_2) : I\to\{x\in\R^2 : \det D_xF_{\rhat}(x)\neq 0\}
  \end{equation*}
  is a solution to the initial value problem
  \begin{subequations}
    \label{eq:h-ode}
    \begin{align}
      \dot{h}(s) &= \big[D_xF_{\rhat}(h(s))\big]\inv\rhat,\label{eq:ode}\\
      h(0)&= x^{(j)}.\label{eq:initial-value}
    \end{align}    
  \end{subequations}
  Then for every $\lambda\in\C\setminus\{0\}$ such that
  $|\lambda|\in I$ the triple $(E(\lambda),N(\lambda),n(\lambda))$
  defined by
  \begin{equation}
    \label{eq:solution-from-IVP}
    E(\lambda)
    := e^{i\theta}\frac{\lambda}{|\lambda|}
    \begin{bmatrix}
      h_1(|\lambda|)\,e^{i\phi_-}\\
      h_2(|\lambda|)\,e^{i\phi_+}
    \end{bmatrix}
    \text{ and }
    \begin{bmatrix}
      N(\lambda)\\n(\lambda)
    \end{bmatrix}
    := y(h(|\lambda|))
  \end{equation}
  is an equilibrium point of system~\eqref{eq:system} with injected
  field $u=\lambda\uhat$.
\end{theorem}

\begin{figure}[t]
  \centering \input{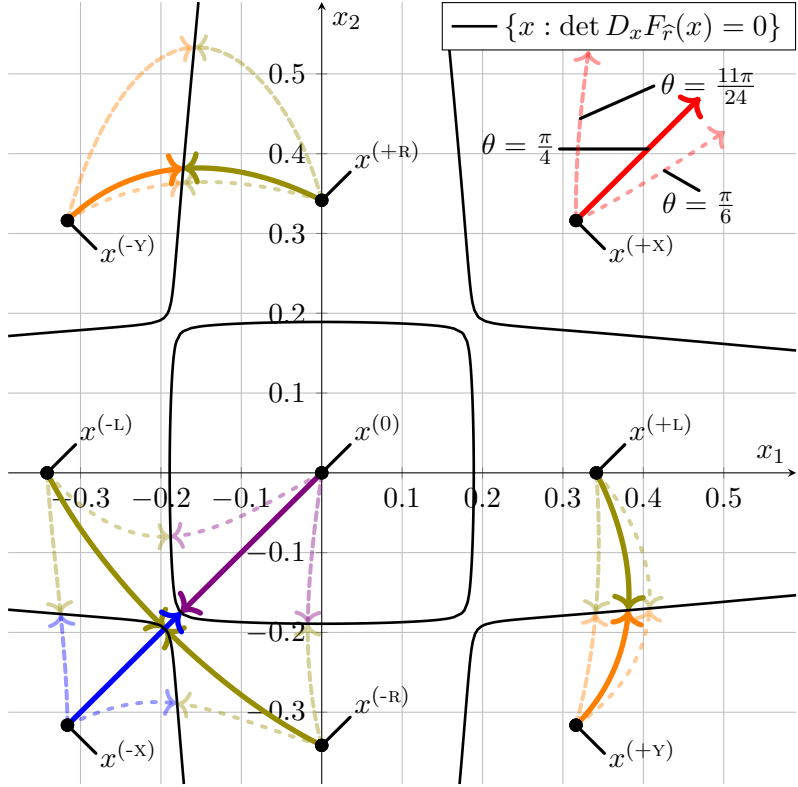}
  
  \caption[Solutions $h$]{Paths traced by solutions
    $h^{(j)}_{\rhat}(s)$, $j\in\J$, of equation~\eqref{eq:F-h-iff} as
    $s\ge 0$ increases. Black dots denote the initial values
    $h^{(j)}_{\rhat}(0)=x^{(j)}$. The paths were solved with
    BifurcationKit.jl~\cite{veltz:hal-02902346}.

    Black lines denote complement of the domain of the initial value
    problem~\eqref{eq:h-ode}. A solution of~\eqref{eq:h-ode} with
    initial value $x^{(j)}$ coincides with $h^{(j)}_{\rhat}(s)$ for as
    long as it does not hit the boundary of the domain (at which point
    the right-hand side of~\eqref{eq:ode} ceases to exist). This means
    that the solution of~\eqref{eq:h-ode} starting from
    $x^{(-\textsc{x})}$ follows the blue path up to the point where
    the path first crosses the black line, and then ends there. All
    other paths can be solved in full from the inital value
    problem~\eqref{eq:h-ode}.

    The solutions $h^{(j)}_{\rhat}$ correspond
    via~\eqref{eq:solution-from-IVP} to the equilibrium points
    (depicted in Figure~\ref{fig:paths}) of a laser with injected
    external optical field. The parameters used are those of
    Figures~\ref{fig:ODE-solution} and~\ref{fig:paths}.}
  \label{fig:h-solution}
\end{figure}

\begin{remark}
  Initial value problem~\eqref{eq:h-ode} is straightforward to solve
  numerically using the explicit expressions for $x^{(j)}$ and
  $D_xF_{\rhat}$ given in~\eqref{eq:xk} and~\eqref{eq:DxF-expression},
  respectively. Therefore Theorem~\ref{thm:solution-from-IVP} provides
  an easy method to trace the trajectories of the equilibrium points
  $E^{(j)}_{\uhat}$ starting from $\lambda=0$ for as long as
  $|\lambda|$ is in the domain $I$ of existence of a solution
  of~\eqref{eq:h-ode}. Also, the asymptotics of $E^{(j)}_{\uhat}$ as
  $\lambda\to 0$ immediately follow from the initial value
  problem~\eqref{eq:h-ode}. On the other hand, if $I$ is a finite
  interval, it may be possible to continue the trajectories even
  beyond the interval $I$. In that case one can use numerical
  continuation techniques, such as pseudo-arclength continuation, to
  solve the functions $h^{(j)}_{\rhat}(s)$ from~\eqref{eq:F-h-iff} and
  use them in~\eqref{eq:solution-from-IVP} instead (cf.\
  Figure~\ref{fig:h-solution}).
\end{remark}
\begin{remark}
  Note that if $\uhat_-=0$, then every $\phi_-\in\R$
  satisfies~\eqref{eq:phi-choice}, and each one of these yields an
  equilibrium point when plugged into~\eqref{eq:solution-from-IVP}. If
  $\uhat_+=0$, an analogous statement holds for $\phi_+$.
\end{remark}

\begin{proof}[Proof of Theorem~\ref{thm:solution-from-IVP}]
  Let $h:I\to\R^2$ solve~\eqref{eq:h-ode}. Then by~\eqref{eq:ode} and
  the chain rule
  \begin{equation}
    \label{eq:F-derivative}
    \frac{d}{ds} F_{\rhat}(s,h(s)) =
    \begin{bmatrix}
      -\rhat & D_xF_{\rhat}(h(s))
    \end{bmatrix}
    \begin{bmatrix}
      1 \\
      \dot{h}(s)
    \end{bmatrix}
    = 0,
  \end{equation}
  so the map $I\ni s\mapsto F_{\rhat}(s,h(s))$ is constant, and
  by~\eqref{eq:initial-value} and Proposition~\ref{prop:zeros-at-0}
  the constant is zero.
  
  Consider $\lambda\neq 0$ such that $|\lambda|\in I$, and choose
  $\phi'_\pm\in\R$ such that
  \begin{equation*}
    e^{i\phi'_\pm} = \frac{\lambda}{|\lambda|}e^{i(\theta+\phi_\pm)}.
  \end{equation*}
  Because $F_{\rhat}(|\lambda|, h(|\lambda|)) = 0$, it follows that
  $x:=h(|\lambda|)$ satisfies $X(y(x))x=|\lambda|\rhat$. Therefore by
  Proposition~\ref{prop:algebraic-version} the triple
  $(E(\lambda),N(\lambda),n(\lambda))$ with
  \begin{equation}
    \label{eq:solution-from-IVP-2}
    E(\lambda) :=
    \begin{bmatrix}
      x_1\, e^{i\phi'_-}\\
      x_2\, e^{i\phi'_+}
    \end{bmatrix}
    \text{ and }
    \begin{bmatrix}
      N(\lambda)\\n(\lambda)
    \end{bmatrix}
    := y(x)
  \end{equation}
  is an equilibrium point of system~\eqref{eq:system} with injected
  field
  \begin{equation}
    \label{eq:lambda-v}
    u := (1+i\alpha)
    |\lambda|
    \begin{bmatrix}
      \rhat_1\,e^{i\phi'_-}\\
      \rhat_2\,e^{i\phi'_+}
    \end{bmatrix}.
  \end{equation}
  Noticing that~\eqref{eq:solution-from-IVP}
  and~\eqref{eq:solution-from-IVP-2} coincide and that the right-hand
  side of~\eqref{eq:lambda-v} is equal to $\lambda\uhat$ finishes the
  proof.
\end{proof}

We are now ready to prove
Theorem~\ref{thm:equilibrium-small-dynamics}:

\begin{proof}[Proof of Theorem~\ref{thm:equilibrium-small-dynamics}]
  We will first prove that there exists a constant $\ell>0$ and nine
  continuous functions $E^{(j)}_{\uhat}$, $j\in\J$, that are of the
  form~\eqref{eq:E-lambda-j}, for which the points $(E,N,n)$ with
  \begin{equation}
    \label{eq:E-N-n}
    E=E^{(j)}_{\uhat}(\lambda)\text{ and }(N,n)=y(|E_-|,|E_+|)
  \end{equation}
  are equilibrium points of system~\eqref{eq:system} with
  $u=\lambda\uhat$, and that satisfy the
  asymptotics~\eqref{eq:E-approximations} as $\lambda\to 0$.
  
  Define
  \begin{equation}
    \label{eq:r-hat-final}
    \rhat := \frac{1}{|1+i\alpha|}
    \begin{bmatrix}
      |\uhat_-|\\
      |\uhat_+|
    \end{bmatrix}
    \in(0,\infty)\times(0,\infty),
  \end{equation}
  and let $\ell>0$ be the constant and
  $h_{\rhat}^{(j)}:(-\ell,\ell)\to\R^2$, $j\in\J$, the smooth
  functions from
  Proposition~\ref{prop:implicit-function-theorem}. Define
  \begin{equation}
    \label{eq:E-lambda}
    E_{\uhat}^{(j)}(\lambda)
    := \frac{\lambda}{|\lambda|}e^{i\theta}
    \begin{bmatrix}
      \frac{\uhat_-}{|\uhat_-|} & 0\\
      0 & \frac{\uhat_+}{|\uhat_+|}
    \end{bmatrix}
    h^{(j)}_{\rhat}(|\lambda|)
    \qquad(j\in\J, 0<|\lambda|<\ell).
  \end{equation}
  Note that if $|\uhat_-|=|\uhat_+|$, then $\rhat_1=\rhat_2$, and for
  $j\in\{0, \pm\textsc{x}\}$ it follows from~\eqref{eq:E-lambda}
  and~\eqref{eq:r1r2} that
  $E_{\uhat}^{(j)}(\lambda) = \rho(\lambda)\uhat$ for some
  $\rho(\lambda)\in\C$.
  
  Fix $j\in\J$. If $0<|\lambda|<\ell$, then
  $x:=h^{(j)}_{\rhat}(|\lambda|)$ satisfies $X(y(x))x=|\lambda|\rhat$,
  and therefore from Proposition~\ref{prop:algebraic-version} it
  follows that a point $(E,N,n)$ defined by~\eqref{eq:E-N-n} is an
  equilibrium point of system~\eqref{eq:system} with
  \begin{equation*}
    u = (1+i\alpha)\frac{\lambda}{|\lambda|}e^{i\theta}
    \begin{bmatrix}
      \frac{\uhat_-}{|\uhat_-|} & 0\\
      0 & \frac{\uhat_+}{|\uhat_+|}
    \end{bmatrix}
    |\lambda|\rhat = \lambda\uhat.
  \end{equation*}
  
  Because the function $h^{(j)}_{\rhat}$ is differentiable, it holds
  that
  \begin{equation}
    \label{eq:h-asymptotic}
    h^{(j)}_{\rhat}(s)
    = h^{(j)}_{\rhat}(0) + s\cdot\frac{d}{ds}h^{(j)}_{\rhat}(0) + o(s)
    \text{ as } s\to 0.
  \end{equation}
  The function $s\mapsto F_{\rhat}(s,h^{(j)}_{\rhat}(s))$ vanishes
  identically, so differentiating it and simplifying
  (see~\eqref{eq:F-derivative}) gives
  \begin{equation*}
    D_xF_{\rhat}(h^{(j)}_{\rhat}(s))\frac{d}{ds}h^{(j)}_{\rhat}(s)
    =\rhat,
  \end{equation*}
  which by Proposition~\ref{prop:F-properties} can be solved at $s=0$
  to yield
  \begin{equation}
    \label{eq:ds}
    \frac{d}{ds}h^{(j)}_{\rhat}(0) = \big[D_xF_{\rhat}(x^{(j)})\big]\inv\rhat.
  \end{equation}

  The matrix $[D_xF_{\rhat}(x^{(j)})]\inv$ in~\eqref{eq:ds} was
  calculated in
  Proposition~\ref{prop:F-properties}. Substituting~\eqref{eq:ds} and
  the value of $h^{(j)}_{\rhat}(0)=x^{(j)}$ from
  Proposition~\ref{prop:zeros-at-0} into~\eqref{eq:h-asymptotic}, and
  then inserting the resulting expression into~\eqref{eq:E-lambda},
  shows that the function $E^{(j)}_{\uhat}$ satisfies
  asymptotics~\eqref{eq:E-approximations} as $\lambda\to 0$. It then
  follows from~\eqref{eq:E-approximations} and the continuity of
  $E_{\uhat}^{(j)}$ that by decreasing $\ell>0$ if necessary, the
  family $\{E_{\uhat}^{(j)}\}_{j\in\J}$ of functions can be made to
  have pairwise distinct values.

  It only remains to prove that if a triple $(E,N,n)$ is an
  equilibrium point of system~\eqref{eq:system} with injected field
  $\lambda\uhat$, where $0<|\lambda|<\ell$, then $E=E^{(j)}(\lambda)$
  for some $j\in\J$, and $(N,n)=y(|E_-|,|E_+|)$. To that end, consider
  an arbitrary equilibrium point $(E,N,n)$ of system~\eqref{eq:system}
  with $u=\lambda\uhat$, where $0<|\lambda|<\ell$. By
  Proposition~\ref{prop:algebraic-version} there exists $x\in\R^2$,
  $r\in[0,\infty)\times[0,\infty)$ and $\phi_\pm\in\R$ such that
  \begin{subequations}
    \begin{align}
      X(y(x))x &= r,\label{eq:necessary-1}\\
      E &=
          \begin{bmatrix}
            x_1\,e^{i\phi_-}\\
            x_2\,e^{i\phi_+}
          \end{bmatrix}
      ,\label{eq:necessary-2}\\
      \begin{bmatrix}
        N\\n
      \end{bmatrix}
               & = y(x),\text{ and}\label{eq:necessary-3}\\
      \lambda\uhat &= (1+i\alpha)
                     \begin{bmatrix}
                       r_1\, e^{i\phi_-}\\
                       r_2\, e^{i\phi_+}
                     \end{bmatrix}
      \label{eq:necessary-4}.
    \end{align}
  \end{subequations}
  
  Equalities~\eqref{eq:necessary-2} and~\eqref{eq:necessary-3} imply
  that $(N,n)=y(|E_-|,|E_+|)$. Also, positivity of the components of
  $r$ together with~\eqref{eq:r-hat-final} and~\eqref{eq:necessary-4}
  imply that $r=|\lambda|\rhat$. Then~\eqref{eq:necessary-1} implies
  that $F_{\rhat}(|\lambda|,x)=0$, so $x=h^{(j)}_{\rhat}(|\lambda|)$
  for some $j\in\J$ by
  Proposition~\ref{prop:implicit-function-theorem}. Finally, dividing
  the components of~\eqref{eq:necessary-4} by their modulus shows that
  \begin{equation*}
    \frac{\lambda}{|\lambda|}\frac{\uhat_\pm}{|\uhat_\pm|}
    = \frac{1+i\alpha}{|1+i\alpha|}e^{i\phi_\pm} = e^{-i\theta}e^{i\phi_\pm}.
  \end{equation*}
  Solving for $e^{i\phi_\pm}$ and inserting these values
  into~\eqref{eq:necessary-2} shows that $E$ is equal to the
  right-hand side of~\eqref{eq:E-lambda}.
\end{proof}

\subsection{Stability of equilibrium points with weak injected fields}
\label{sec:stability}

\begin{figure}[tp]
  \centering \input{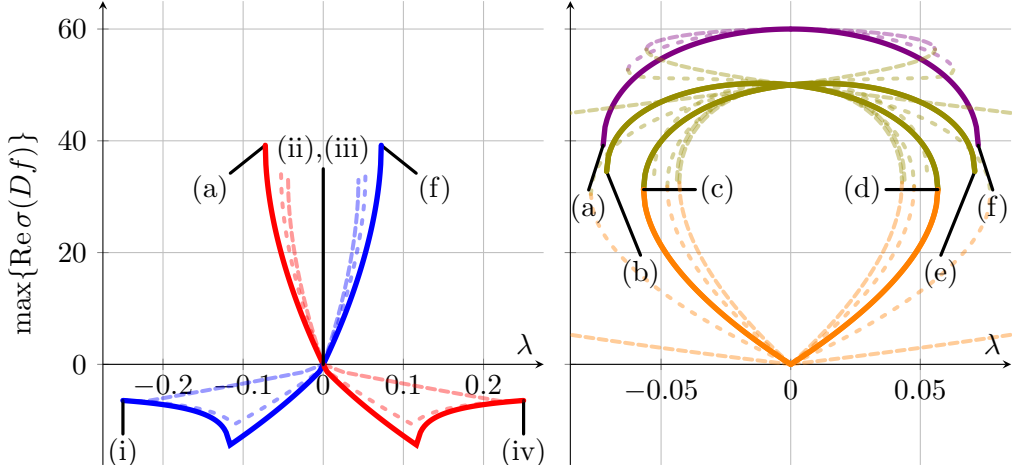}
  \caption[Stability]{Linear stability analysis of equilibrium points
    of a laser subject to external optical injection
    $u=\lambda\uhat$. Each line represents an equilibrium point, and
    $\max\{\re\sigma(Df)\}$ denotes the maximum real part of
    eigenvalues of the linearized system at an equilibrium point. A
    positive value indicates that the equilibrium point is unstable,
    while a negative value indicates that the equilibrium point is
    asymptotically stable.  The parameters used, the color and style
    of the lines, as well as the labels {(a)}--{(f)} and {(i)}--{(iv)}
    match those of Figure~\ref{fig:paths}.
    
    At $\lambda=0$ the blue line from {(i)} to {(f)} and the red line
    from {(a)} to {(iv)} change signs, this corresponds to a jump of
    the stable equilibrium point from {(ii)} to {(iii)} in
    Figure~\ref{fig:paths}. The other lines with shorter intervals of
    existence are positive for all $\lambda\neq 0$, they are displayed
    on the axis on the right-hand side (note the different scales on
    the $\lambda$-axes).}
  \label{fig:stability}
\end{figure}
  
In this section, we consider stability properties of the nine
equilibrium points from
Theorem~\ref{thm:equilibrium-small-dynamics}. We will prove in
Theorem~\ref{thm:stability} below that if $\alpha=0$ and the injected
field $u$ in system~\eqref{eq:system} is sufficiently weak, then the
system has exactly one \emph{asymptotically stable} (in the sense of
Lyapunov) equilibrium point, while the remaining equilibrium points
are \emph{unstable} (for the definitions of asymptotic stability and
instability of an equilibrium point, we refer the reader
to~\cite{MR1071170}).

By splitting the complex-valued functions $E_\pm(t)$ into their real
and imaginary parts, i.e., writing
$E_\pm(t)= E_\pm^{(\mathrm{re})}(t) + iE_\pm^{(\mathrm{im})}(t)$ with
${E_\pm^{(\mathrm{re})}(t), E_\pm^{(\mathrm{im})}(t)}\in\R$, we can
write system~\eqref{eq:system} in terms of real-valued functions as
\begin{equation}
  \label{eq:real-system}
  \frac{d}{dt}(E_-^{(\mathrm{re})},E_+^{(\mathrm{re})},E_-^{(\mathrm{im})},E_+^{(\mathrm{im})},N,n)
  = f(E_-^{(\mathrm{re})},E_+^{(\mathrm{re})},E_-^{(\mathrm{im})},E_+^{(\mathrm{im})},N,n),
\end{equation}
where the function $f:\R^6\to\R^6$ is determined by
system~\eqref{eq:system}. A calculation shows that $Df$, the Jacobian
matrix of $f$, is given by the block matrix
\begin{equation}
  \begin{split}
    &Df(E_-^{(\mathrm{re})},E_+^{(\mathrm{re})},E_-^{(\mathrm{im})},E_+^{(\mathrm{im})},N,n) =\\
    &-\!\!
    \begin{bmatrix}
      \label{eq:Df}
      \kappa X(N,n)&-\alpha\kappa X(N,n)&-\kappa(F^{(\mathrm{re})}-\alpha F^{(\mathrm{im})})\\
      \alpha\kappa X(N,n)&\kappa X(N,n)&-\kappa (\alpha F^{(\mathrm{re})}+F^{(\mathrm{im})})\\
      2\gamma(F^{(\mathrm{re})})^T(I_2-X(N,n))&2\gamma(F^{(\mathrm{im})})^T(I_2-X(N,n))&\gamma
      Y(E)
    \end{bmatrix}\!\!,
  \end{split}
\end{equation}
where the superscript $T$ denotes the transpose of a matrix, and
\begin{equation*}
  F^{(j)} :=
  \begin{bmatrix*}[r]
    E_-^{(j)} & -E_-^{(j)}\\
    E_+^{(j)} & E_+^{(j)}
  \end{bmatrix*}
  \qquad(j\in\{{\mathrm{im},\mathrm{re}}\}).
\end{equation*}

We proved in Theorem~\ref{thm:solution-from-IVP} a method for
calculating numerical values for the nine equilibrium points of
system~\eqref{eq:system} from
Theorem~\ref{thm:equilibrium-small-dynamics}. Inserting the value of
an equilibrium point into the expression~\eqref{eq:Df} for $Df$ and
finding the eigenvalues of the so obtained $6\times 6$-matrix is an
easy numerical method to test the stability of the equilibrium
point. Recall that if all the eigenvalues of $Df$ at an equilibrium
point have strictly negative real parts, then the equilibrium point is
asymptotically stable, while if at least one of the eigenvalues has a
strictly positive real part, then the equilibrium point is
unstable~\cite{MR1071170}. Only if none of the eigenvalues have
strictly positive real parts but at least one of them has real part
equal to zero, then this test for stability is inconclusive.

In Figure~\ref{fig:stability} we have used above test to determine
stability of the equilibrium points on Figure~\ref{fig:paths}. As
illustrated in Figure~\ref{fig:stability}, of the nine equilibrium
points depicted in Figure~\ref{fig:paths} that correspond to an
injected field $u=\lambda\widehat{u}$, for each
$\lambda\in[-1/4,1/4]\setminus{\{0\}}$ and $\widehat{u}$ exactly one
of the points is asymptotically stable, while the others are unstable.

\begin{lemma}
  \label{lemma:eigs-of-Df}
  Assume that $\alpha=0$ in system~\eqref{eq:system}, and consider the
  Jacobian matrix $Df$ of the corresponding
  system~\eqref{eq:real-system}. (An expression for $Df$ is given
  at~\eqref{eq:Df}.)
  \begin{enumerate}[(i)]
  \item For arbitrary numbers
    ${E_\pm^{(\mathrm{re})}, E_\pm^{(\mathrm{im})}, N, n}\in\R$, and
    for the matrix
    \begin{equation*}
      Df = Df(E_-^{(\mathrm{re})},E_+^{(\mathrm{re})},E_-^{(\mathrm{im})},E_+^{(\mathrm{im})},N,n)
    \end{equation*}
    the following hold:
    \begin{align}
      (E_-^{(\mathrm{im})}, 0, -E_-^{(\mathrm{re})}, 0, 0, 0)&\in\ker(Df+\kappa [1-(N-n)]I_6)\text{ and}
                                                               \label{eq:eigs-of-Df-1}\\
      (0, E_+^{(\mathrm{im})}, 0, -E_+^{(\mathrm{re})}, 0, 0)&\in\ker(Df+\kappa [1-(N+n)]I_6).
                                                               \label{eq:eigs-of-Df-2}
    \end{align}
  \item Let $\theta_1$ and $\theta_2$ be the two roots of the
    polynomial
    \begin{equation}
      \label{eq:eigs-of-Df-p1}
      s^2 + \gamma\mu s + 2\kappa\gamma(\mu-1),
    \end{equation}
    and let $\theta_3$ and $\theta_4$ be the two roots of the
    polynomial
    \begin{equation}
      \label{eq:eigs-of-Df-p2}
      s^2 + \gamma(\delta+\mu-1)s + 2\kappa\gamma(\mu-1).
    \end{equation}
    Furthermore, let
    ${E_\pm^{(\mathrm{re})},E_\pm^{(\mathrm{im})}}\in\R$ be any
    numbers such that
    \begin{equation*}
      |E_-^{(\mathrm{re})}+iE_-^{(\mathrm{im})}|=|E_+^{(\mathrm{re})}+iE_+^{(\mathrm{im})}|=\sqrt{
        \frac{\mu-1}{2}}.
    \end{equation*}
    Then for the matrix
    \begin{equation*}
      Df = Df(E_-^{(\mathrm{re})},E_+^{(\mathrm{re})},E_-^{(\mathrm{im})},E_+^{(\mathrm{im})},1,0)
    \end{equation*}
    the following holds:
    \begin{align}
      (E_-^{(\mathrm{re})}, E_+^{(\mathrm{re})}, E_-^{(\mathrm{im})},  E_+^{(\mathrm{im})}, \theta_j/\kappa, 0)
      &\in\ker(Df-\theta_jI_6), \, j = 1, 2,\text{ and}
        \label{eq:eigs-of-Df-3}\\
      (-E_-^{(\mathrm{re})}, E_+^{(\mathrm{re})}, -E_-^{(\mathrm{im})}, E_+^{(\mathrm{im})}, 0, \theta_j/\kappa)
      &\in\ker(Df-\theta_jI_6), \, j = 3, 4.
        \label{eq:eigs-of-Df-4}
    \end{align}
  \end{enumerate}
\end{lemma}
\begin{proof}
  The straightforward calculation using expression~\eqref{eq:Df} for
  $Df$ is omitted.
\end{proof}

Given $\uhat\in\C^2$ such that $\uhat_-\neq 0$ and $\uhat_+\neq 0$,
let $E_{\uhat}^{(j)}$, $j\in\J$, be the functions from
Theorem~\ref{thm:equilibrium-small-dynamics}. By {(ii)} of
Proposition~\ref{prop:algebraic-version}, if
$(E,N,n)\in\C^2\times\R\times\R$ is an equilibrium point of
system~\eqref{eq:system}, then $(N,n)=y(|E_-|,|E_+|)$. Therefore we
can define functions
\begin{equation*}
  \lambda\mapsto N_{\uhat}^{(j)}(\lambda),\,
  \lambda\mapsto n_{\uhat}^{(j)}(\lambda),\text{ and }
  \lambda\mapsto(Df)_{\uhat}^{(j)}(\lambda)
\end{equation*}
in a punctured neighborhood of the origin of the complex plane by
requiring that the point
$(E_{\uhat}^{(j)}(\lambda),N_{\uhat}^{(j)}(\lambda),n_{\uhat}^{(j)}(\lambda))\in\C^2\times\R\times\R$
is an equilibrium point of system~\eqref{eq:system}, and that
$(Df)_{\uhat}^{(j)}(\lambda)$ is the Jacobian matrix of
system~\eqref{eq:real-system} at that point. In other words, if
$\lambda\neq 0$ is sufficiently small and
${E_\pm^{(\mathrm{re})},E_\pm^{(\mathrm{im})}}\in\R$ are such that
$E_{\uhat}^{(j)}(\lambda)=(E_-^{(\mathrm{re})}+iE_-^{(\mathrm{im})},E_+^{(\mathrm{re})}+iE_+^{(\mathrm{im})})$,
then
\begin{align}
  \begin{bmatrix*}
    N_{\uhat}^{(j)}(\lambda)\\
    n_{\uhat}^{(j)}(\lambda)
  \end{bmatrix*}
  &= y(|E_-^{(\mathrm{re})}+iE_-^{(\mathrm{im})}|,|E_+^{(\mathrm{re})}+iE_+^{(\mathrm{im})}|),\text{ and}\label{eq:Nn-def}\\
  (Df)_{\uhat}^{(j)}(\lambda)
  &= Df(E_-^{(\mathrm{re})},E_+^{(\mathrm{re})},E_-^{(\mathrm{im})},E_+^{(\mathrm{im})},N_{\uhat}^{(j)}(\lambda),n_{\uhat}^{(j)}(\lambda)).\label{eq:Df-lambda}
\end{align}
We call an equilibrium point
$(E_{\uhat}^{(j)}(\lambda),N_{\uhat}^{(j)}(\lambda),n_{\uhat}^{(j)}(\lambda))$
the \emph{equilibrium point corresponding to
  $E_{\uhat}^{(j)}(\lambda)$}.

We can now prove instability for five of the equilibrium points from
Theorem~\ref{thm:equilibrium-small-dynamics}:

\begin{lemma}
  \label{lemma:unstable-points}
  Assume $\alpha=0$ in system~\eqref{eq:system}. Fix $\uhat\in\C^2$
  with $\uhat_-\neq 0$ and $\uhat_+\neq 0$, and let $\ell>0$ and
  $E_{\uhat}^{(j)}$, $j\in\J$, be as in
  Theorem~\ref{thm:equilibrium-small-dynamics}. Then there exists
  $0<\ell_0\le\ell$ such that if $0<|\lambda|<\ell_0$, then the
  equilibrium points corresponding to $E_{\uhat}^{(j)}(\lambda)$ with
  $j\in\{\textsc{0},\pm\textsc{l},\pm\textsc{r}\}$ are unstable.
\end{lemma}
\begin{proof}
  Choose $j\in\J$ and sufficiently small $\lambda\neq 0$, and set
  $(E_-,E_+):=E_{\uhat}^{(j)}(\lambda)$.  By {(i)} of
  Lemma~\ref{lemma:eigs-of-Df}, the number $-\kappa[1-(N-n)]$ is an
  eigenvalue of $(Df)_{\uhat}^{(j)}(\lambda)$ if $E_-\neq 0$, and
  $-\kappa[1-(N+n)]$ is an eigenvalue of $(Df)_{\uhat}^{(j)}(\lambda)$
  if $E_+\neq 0$. It follows from the
  asymptotics~\eqref{eq:E-approximations} that there exists
  $0<\ell_1\le\ell$ such that $E_-\neq 0$ and $E_+\neq 0$ if
  $0<|\lambda|<\ell_1$, and therefore the numbers
  $-\kappa[1-(N\pm n)]$ are eigenvalues of
  $(Df)_{\uhat}^{(j)}(\lambda)$ for $0<|\lambda|<\ell_1$.
  
  The limits of $N_{\uhat}^{(j)}(\lambda)$ and
  $n_{\uhat}^{(j)}(\lambda)$ as $\lambda\to 0$ can be calculated using
  the asymptotics~\eqref{eq:E-approximations} of
  $E_{\uhat}^{(j)}(\lambda)$ and~\eqref{eq:Nn-def}. In particular,
  \begin{align*}
    \lim_{\lambda\to 0} -\kappa\big[1-(N_{\uhat}^{(\textsc{0})}(\lambda)\pm n_{\uhat}^{(\textsc{0})}(\lambda))\big]
    &= \kappa(\mu-1) > 0,\\
    \lim_{\lambda\to 0} -\kappa\big[1-(N_{\uhat}^{(\pm\textsc{l})}(\lambda)+ n_{\uhat}^{(\pm\textsc{l})}(\lambda))\big]
    &= 2\kappa\,\frac{\mu-1}{1+\delta}>0,\\
    \lim_{\lambda\to 0} -\kappa\big[1-(N_{\uhat}^{(\pm\textsc{r})}(\lambda)-n_{\uhat}^{(\pm\textsc{r})}(\lambda))\big] &= 2\kappa\,\frac{\mu-1}{1+\delta}>0.
  \end{align*}
  It follows that there exists $0<\ell_0\le\ell_1$ such that if
  $0<|\lambda|<\ell_0$ and
  $j\in\{\textsc{0},\pm\textsc{l},\pm\textsc{r}\}$, then at least one
  of the eigenvalues of $(Df)_{\uhat}^{(j)}(\lambda)$ is strictly
  positive.

  We have shown that the linearization $(Df)_{\uhat}^{(j)}(\lambda)$
  of system~\eqref{eq:real-system} at an equilibrium point
  corresponding to $E_{\uhat}^{(j)}(\lambda)$ with
  $0<|\lambda|<\ell_0$ and
  $j\in\{\textsc{0},\pm\textsc{l},\pm\textsc{r}\}$ has at least one
  strictly positive eigenvalue. Therefore the nonlinear
  system~\eqref{eq:system} is unstable at such a
  point~\cite[Theorem~{15.6}]{MR1071170}.
\end{proof}

Let $\C^6_{\mathrm{sym}}$ denote the quotient space of $\C^6$ by the
equivalence relation that identifies vectors whose coordinates are
permutations of each other, and let
\begin{equation}
  \label{eq:map-sigma}
  \sigma:\C^{6\times 6}\to\C^6_{\mathrm{sym}} 
\end{equation}
denote the map that takes a matrix to the unordered $6$-tuple of its
eigenvalues (repeated according to their algebraic multiplicities).
Then $(\C^6_{\mathrm{sym}},d)$ is a metric space with the
\emph{optimal matching distance}~\cite{MR1477662}
\begin{equation*}
  d([a],[b]) := \min_\beta\max_{1\le k\le 6} |a_k-b_{\beta(k)}|,
\end{equation*}
where $[a]$ and $[b]$ denote the equivalence classes of ${a,b}\in\C^6$
in $\C^6_{\mathrm{sym}}$, and the minimum is taken over all
permutations $\beta$ of $\{1,2,\ldots,6\}$. The map $\sigma$ is
continuous in this topology~\cite{MR1477662}.

Let $\uhat=(\uhat_-,\uhat_+)\in\C^2$ be such that $\uhat_-\neq 0$ and
$\uhat_+\neq 0$. For $\lambda\neq 0$ and
$j\in\{\pm\textsc{x}, \pm\textsc{y}\}$ define
\begin{equation}
  \label{eq:H-def}
  H_{\uhat}^{(j)}(\lambda) :=
  Df(E_-^{(\mathrm{re})},E_+^{(\mathrm{re})},E_-^{(\mathrm{im})},E_+^{(\mathrm{im})},1,0),
\end{equation}
where the arguments $E_\pm^{(\mathrm{re})}\in\R$ and
$E_\pm^{(\mathrm{im})}\in\R$ are defined by
\begin{equation*}
  \begin{bmatrix}
    E_-^{(\mathrm{re})}+iE_-^{(\mathrm{im})}\\
    E_+^{(\mathrm{re})}+iE_+^{(\mathrm{im})}
  \end{bmatrix}
  :=
  \begin{cases}
    \pm\frac{\lambda}{|\lambda|}\sqrt{\frac{\mu-1}{2}}
    \begin{bmatrix}
      \uhat_-/|\uhat_-|\\
      \uhat_+/|\uhat_+|
    \end{bmatrix},\text{ if } j = \pm\textsc{x},\\[1em]
    \pm\frac{\lambda}{|\lambda|}\sqrt{\frac{\mu-1}{2}}
    \begin{bmatrix*}[r]
      \uhat_-/|\uhat_-|\\
      -\uhat_+/|\uhat_+|
    \end{bmatrix*},\text{ if } j = \pm\textsc{y}.
  \end{cases}
\end{equation*}
In other words, $H_{\uhat}^{(j)}(\lambda)$ is defined as
$(Df)_{\uhat}^{(j)}(\lambda)$ in~\eqref{eq:Df-lambda}, except that
$E_{\uhat}^{(j)}(\lambda)$, $N_{\uhat}^{(j)}(\lambda)$ and
$n_{\uhat}^{(j)}(\lambda)$ are replaced by their zeroth order
approximations from~\eqref{eq:E-approximations} (as we are considering
the case $\alpha=0$, we have $e^{i\theta}=1$
in~\eqref{eq:E-approximations}).

Our plan is to determine stability of the remaining equilibrium points
corresponding to $E_{\uhat}^{(j)}(\lambda)$ with
$j\in\{\pm\textsc{x}, \pm\textsc{y}\}$ by finding all eigenvalues of
$(Df)_{\uhat}^{(j)}(\lambda)$. In the following lemma we will first
show that for small $\lambda\neq0$ the eigenvalues of
$H_{\uhat}^{(j)}(\lambda)$ approximate those of
$(Df)_{\uhat}^{(j)}(\lambda)$, and after that in
Lemma~\ref{lemma:eigs-of-H} we will determine the eigenvalues of
$H_{\uhat}^{(j)}(\lambda)$. Combining these results will then make it
possible for us to conclude stability of the equilibrium points.

\begin{lemma}
  \label{lemma:Df-limit}
  For every $j\in\{\pm\textsc{x},\pm\textsc{y}\}$,
  \begin{equation}
    \label{eq:limit-d}
    \lim_{\lambda\to 0} d\big(\sigma\big((Df)_{\uhat}^{(j)}(\lambda)\big),\sigma\big(H_{\uhat}^{(j)}(\lambda)\big)\big) = 0.
  \end{equation}
  Here $d$ is the optimal matching distance on $\C^6_{\mathrm{sym}}$
  and $\sigma$ is the map~\eqref{eq:map-sigma}.
\end{lemma}
\begin{proof}
  A calculation shows that for every
  $j\in\{\pm\textsc{x},\pm\textsc{y}\}$,
  \begin{equation}
    \label{eq:limit-GH}
    \lim_{\lambda\to 0}
    \big\|(Df)_{\uhat}^{(j)}(\lambda)-H_{\uhat}^{(j)}(\lambda)\big\|
    = 0.
  \end{equation}
  There exists numbers $r>0$ and $R>0$ such that if $0<|\lambda|<r$
  and $j\in\{\pm\textsc{x},\pm\textsc{y}\}$, then
  $(Df)_{\uhat}^{(j)}(\lambda)\in\overbar{B}_R$ and
  $H_{\uhat}^{(j)}(\lambda)\in\overbar{B}_R$, where
  $\overbar{B}_R\subset\R^{6\times 6}$ is the closed ball of radius
  $R$ centered at the origin. Because the continuous map $\sigma$ is
  uniformly continuous on the compact set $\overbar{B}_R$,
  from~\eqref{eq:limit-GH} it follows that the
  limit~\eqref{eq:limit-d} holds.
\end{proof}

\begin{lemma}
  \label{lemma:eigs-of-H}
  Let $\theta_j\in\C$, $j\in\{1,2,3,4\}$, be the roots in {(ii)} of
  Lemma~\ref{lemma:eigs-of-Df}. If
  $j\in\{\pm\textsc{x},\pm\textsc{y}\}$ and $\lambda\neq 0$, then
  $(0,0,\theta_1,\theta_2,\theta_3,\theta_4)$ is a sequence of all
  eigenvalues of $H_{\uhat}^{(j)}(\lambda)$ (repeated according to
  their algebraic multiplicities).
\end{lemma}
\begin{proof}
  Because $N=1$ and $n=0$ in the definition~\eqref{eq:H-def} of
  $H_{\uhat}^{(j)}(\lambda)$, {(i)} of Lemma~\ref{lemma:eigs-of-Df}
  implies that zero is an eigenvalue of $H_{\uhat}^{(j)}(\lambda)$. By
  {(ii)} of the same lemma, also the four roots $\theta_j$ are
  eigenvalues of $H_{\uhat}^{(j)}(\lambda)$.

  If $\theta_1\neq\theta_2$ and $\theta_3\neq\theta_4$, it can be
  calculated that the six vectors on the left-hand sides
  of~\eqref{eq:eigs-of-Df-1}, \eqref{eq:eigs-of-Df-2},
  \eqref{eq:eigs-of-Df-3}, and \eqref{eq:eigs-of-Df-4} form a linearly
  independent set. It follows that in this case
  $(0,0,\theta_1,\theta_2,\theta_3,\theta_4)$ is a sequence of all
  eigenvalues of $H_{\uhat}^{(j)}(\lambda)$ (repeated according to
  their algebraic multiplicities).

  If $\theta_1=\theta_2$ or $\theta_3=\theta_4$ we proceed as
  follows. So far $\gamma>0$ has been fixed, let us now temporarily
  write $H_{\uhat}^{(j)}(\lambda,\gamma)$ to consider
  $H_{\uhat}^{(j)}$ as a function of both $\lambda$ and
  $\gamma>0$. Also, denote by $\theta_j(\gamma)$ the roots of the
  polynomials~\eqref{eq:eigs-of-Df-p1} and~\eqref{eq:eigs-of-Df-p2}
  for given $\gamma$ (in arbitrary order).
  
  For $\lambda\neq 0$ fixed, both of the maps
  $(0,\infty)\ni\gamma\mapsto\sigma(H_{\uhat}^{(j)}(\lambda,\gamma))\in\C^6_{\mathrm{sym}}$
  and
  $(0,\infty)\ni\gamma\mapsto[(0, 0, \theta_1(\gamma),
  \theta_2(\gamma),
  \theta_3(\gamma),\theta_4(\gamma))]\in\C^6_{\mathrm{sym}}$ are
  continuous. By the first part of the proof these maps agree except
  possibly for the finite set of $\gamma$ where one of the
  polynomials~\eqref{eq:eigs-of-Df-p1} and~\eqref{eq:eigs-of-Df-p2}
  has a double root. But by continuity they then agree everywhere.
\end{proof}

We can now prove the main result of this section.

\begin{theorem}
  \label{thm:stability}
  Consider system~\eqref{eq:system} under the assumption that
  $\alpha=0$ and that the injected field $u$ is of the form
  $u=\lambda\uhat$, where $\lambda\in\C\setminus\{0\}$ and
  $\uhat=(\uhat_-,\uhat_+)\in\C^2$ satisfies $\uhat_-\neq 0$ and
  $\uhat_+\neq 0$. With reference to
  Theorem~\ref{thm:equilibrium-small-dynamics}, let $\ell>0$ be a
  constant and $E_{\uhat}^{(j)}(\lambda)$, $j\in\J$, the functions
  with asymptotics~\eqref{eq:E-approximations} such that for
  $0<|\lambda|<\ell$ they determine the nine equilibrium points of
  system~\eqref{eq:system} with injected field $u=\lambda\uhat$.

  There exists a constant $0<\ell_0\le\ell$ such that for every
  $0<|\lambda|<\ell_0$ the equilibrium point corresponding to
  $E_{\uhat}^{(+\textsc{x})}(\lambda)$ is asymptotically stable, and
  the other eight equilibrium points corresponding to
  $E_{\uhat}^{(j)}(\lambda)$ with
  $j\in\{\textsc{0},\pm\textsc{l},\pm\textsc{r},-\textsc{x},\pm\textsc{y}\}$
  are unstable.
\end{theorem}
\begin{proof}
  By Lemma~\ref{lemma:unstable-points} we know that the equilibrium
  points corresponding to $E_{\uhat}^{(j)}(\lambda)$ with
  $j\in\{\textsc{0},\pm\textsc{l},\pm\textsc{r}\}$ and $\lambda\neq0$
  sufficiently small are unstable. By decreasing $\ell>0$ if
  necessary, we can assume that this is the case for all
  $0<|\lambda|<\ell$.

  To prove the theorem, we will show that for sufficiently small
  $\lambda\neq0$ all of the eigenvalues of
  $(Df)_{\uhat}^{(+\textsc{x})}(\lambda)$ have strictly negative real
  parts, and that at least one of the eigenvalues of each of
  $(Df)_{\uhat}^{(j)}(\lambda)$ with
  $j\in\{-\textsc{x},\pm\textsc{y}\}$ has a strictly positive real
  part. By~\cite[Theorem~{15.6}]{MR1071170} this will imply the
  result.

  Let $\theta_i$, $i\in\{1,2,3,4\}$, be the roots of the
  polynomials~\eqref{eq:eigs-of-Df-p1} and~\eqref{eq:eigs-of-Df-p2} in
  Lemma~\ref{lemma:eigs-of-Df}. Because all of the coefficients in the
  polynomials are strictly positive, $\re\theta_i<0$ for every
  $i$. Therefore it is possible to find a radius $r>0$ such that
  $\cup_{i=1}^4B_r(\theta_i)\subset\C_-:=\{z\in\C:\re z<0\}$, and such
  that this union is disjoint from $B_r(0)$. Here $B_r(z)\subset\C$
  denotes the open disk of radius $r$ centered at $z\in\C$.

  Fix $j\in\{\pm\textsc{x},\pm\textsc{y}\}$. By
  Lemmas~\ref{lemma:Df-limit} and~\ref{lemma:eigs-of-H} and the
  definition of the optimal matching distance $d$, we can find
  $0<\ell_1\le\ell$ such that if $0<|\lambda|<\ell_1$, then
  $(Df)_{\uhat}^{(j)}(\lambda)$ has two eigenvalues in $B_r(0)$ and
  four eigenvalues in $\cup_{i=1}^4B_r(\theta_i)$. A calculation shows
  that
  \begin{equation*}
    \lim_{\lambda\to 0} \kappa\big[1-(N_{\uhat}^{(j)}(\lambda)\pm n_{\uhat}^{(j)}(\lambda))\big] = 0,
  \end{equation*}
  so by {(i)} of Lemma~\ref{lemma:eigs-of-Df} the two eigenvalues of
  $(Df)_{\uhat}^{(j)}(\lambda)$ contained in $B_r(0)$ are
  \begin{equation}
    \label{eq:relevant-eigs}
    -\kappa\big[1-(N_{\uhat}^{(j)}(\lambda)\pm n_{\uhat}^{(j)}(\lambda))\big].
  \end{equation}
  Because $\cup_{i=1}^4B_r(\theta_i)\subset\C_-$, only the two
  eigenvalues~\eqref{eq:relevant-eigs} are relevant for determining
  the stability for small $\lambda\neq0$.

  Consider Theorem~\ref{thm:solution-from-IVP} and let $h$ be a
  solution to the initial value problem~\eqref{eq:h-ode}. For
  sufficiently small $\lambda\neq0$ let $E(\lambda)$, $N(\lambda)$ and
  $n(\lambda)$ be defined in terms of $h$
  by~\eqref{eq:solution-from-IVP}. Then by
  Theorem~\ref{thm:equilibrium-small-dynamics} the vector $E(\lambda)$
  is equal to $E_{\uhat}^{(k)}(\lambda)$ for some $k\in\J$, and an
  inspection shows that $k=j$ is the only possibility. If $y_1$ and
  $y_2$ are the component functions of the function $y$
  from~\eqref{eq:def-y}, i.e., $y(x)=(y_1(x),y_2(x))$, above implies
  that
  \begin{equation}
    \label{eq:eig-equality}
    -\kappa\big[1-(N_{\uhat}^{(j)}(\lambda)\pm n_{\uhat}^{(j)}(\lambda))\big]
    = -\kappa\big[1-(y_1(h(|\lambda|))\pm y_2(h(|\lambda|)))\big].
  \end{equation}

  The functions $s\mapsto y_k\circ h(s)$ are defined and
  differentiable in a neighborhood of the origin, and
  \begin{equation}
    \label{eq:dds}
    \begin{split}
      \frac{d}{ds}\big(y_1\circ h\pm y_2\circ h\big)(0)
      &=\nabla (y_1\pm y_2)(h(0))\cdot \frac{d}{ds} h(0)\\
      &=\nabla (y_1\pm
      y_2)(x^{(j)})\cdot\big[D_xF_{\rhat}(x^{(j)})\big]\inv\rhat,
    \end{split}
  \end{equation}
  where $\rhat = (|\uhat_-|, |\uhat_+|)\in(0,\infty)\times(0,\infty)$.
  Calculating the gradient and applying the value of
  $[D_xF_{\rhat}(x^{(j)})]\inv$ obtained in {(i)} of
  Proposition~\ref{prop:F-properties} to~\eqref{eq:dds}, we can
  calculate that
  \begin{align}
    \frac{d}{ds}\big(-\kappa\big[1-(y_1\circ h-y_2\circ h)\big]\big)(0)
    &=-\left(\frac{2\kappa|\uhat_-|}{\mu-1}\right)x_1^{(j)},\text{ and}\label{eq:d-eig-1}\\
    \frac{d}{ds}\big(-\kappa\big[1-(y_1\circ h+y_2\circ h)\big]\big)(0)
    &=-\left(\frac{2\kappa|\uhat_+|}{\mu-1}\right)x_2^{(j)}\label{eq:d-eig-2}.
  \end{align}
      
  The numbers in the parenthesis on the right-hand sides
  of~\eqref{eq:d-eig-1} and~\eqref{eq:d-eig-2} are nonzero and
  positive. If $j=+\textsc{x}$, then $x_1^{(j)}>0$ and $x_2^{(j)}>0$,
  so both~\eqref{eq:d-eig-1} and~\eqref{eq:d-eig-2} are strictly
  negative. This and~\eqref{eq:eig-equality} imply that there exists
  $0<\ell_0\le\ell_1$ such that for $0<|\lambda|<\ell_0$,
  \begin{equation*}
    -\kappa\big[1-(N_{\uhat}^{(+\textsc{x})}(\lambda)\pm n_{\uhat}^{(+\textsc{x})}(\lambda))\big]<0.
  \end{equation*}
  Therefore for these $\lambda$ these two eigenvalues of
  $(Df)_{\uhat}^{(+\textsc{x})}(\lambda)$ are strictly negative, and
  consequently the equilibrium point corresponding to
  $E_{\uhat}^{(+\textsc{x})}(\lambda)$ is asymptotically stable.

  If $j\in\{-\textsc{x},\pm\textsc{y}\}$, then at least one of the
  nonzero numbers $x_1^{(j)}$ and $x_2^{(j)}$ in~\eqref{eq:d-eig-1}
  and~\eqref{eq:d-eig-2} is negative. An analogous reasoning as above
  shows that by decreasing $\ell_0>0$ if necessary, we can conclude
  that for $0<|\lambda|<\ell_0$ at least one of the
  eigenvalues~\eqref{eq:relevant-eigs} of
  $(Df)_{\uhat}^{(j)}(\lambda)$ is strictly positive, and therefore
  the equilibrium point corresponding to $E_{\uhat}^{(j)}(\lambda)$ is
  unstable.
\end{proof}

\subsection{Equilibrium points with strong injected fields}
\label{sec:strong-fields}

In this section, we consider equilibrium points of
system~\eqref{eq:system} under the assumption that the injected
electric field $u$ is strong (large in magnitude). We assume that the
injected field is of the form
\begin{equation*}
  u=\lambda\uhat,
\end{equation*}
where $\lambda\in\C$ is a large parameter and
$\uhat=(\uhat_-,\uhat_+)\in\C^2$ satisfies $\uhat_-\neq 0$ and
$\uhat_+\neq 0$, and we are interested in the behavior of the
equilibrium points as a function of the parameter $\lambda$.

For a number $0<\eta<1$ and a vector $\rhat\in\R^2$ such that
\begin{equation}
  \label{eq:rhat-assumptions}
  \rhat_1>0,\, \rhat_2>0,\text{ and } |\rhat|=1,
\end{equation}
let us define the compact set
\begin{equation*}
  K(\eta,\rhat) := \left\{
    w\in\R^2 : w\cdot \rhat \ge \eta|w|\text{ and } \frac{1}{2}\le|w|\le \frac{3}{2}
  \right\}.
\end{equation*}
We will prove that given the vector $\rhat\in\R^2$, we can choose a
number $\eta=\eta(\rhat)\in(0,1)$ and a constant $L=L(\rhat)>0$ so
that for the function $F_{\rhat}$ defined in~\eqref{eq:def-F} the
following holds: If $s\ge L$, then
\begin{enumerate}[(i)]
\item $F_{\rhat}(s,x)=0$ implies $x\in sK(\eta,\rhat)$, and
\item the map $\R^2\ni x\mapsto x-F_{\rhat}(s,x)\in\R^2$ maps
  $sK(\eta,\rhat)$ contractively into itself.
\end{enumerate}

Recall that by Proposition~\ref{prop:algebraic-version} the zeros of
$F_{\rhat}(s,\cdot)$ and the equilibrium points of
system~\eqref{eq:system} are in one-to-one correspondence. Once {(i)}
and {(ii)} are proved, we can conclude from {(i)} that for $s\ge L$
every zero of $F_{\rhat}(s,\cdot)$ is contained in $sK(\eta,\rhat)$,
and from {(ii)} and the Banach fixed-point theorem that there exists
exactly one such zero in $sK(\eta,\rhat)$. From this it follows that
if the injected field $u$ is strong enough, then there exists a unique
equilibrium point of system~\eqref{eq:system}.

\begin{lemma}
  \label{lemma:i-of-large-lemma}
  Let $0<\eta<1$. There exists a constant $L=L(\eta)>0$ such that if
  $s\ge L$, $\rhat\in\R^2$ is a vector that
  satisfies~\eqref{eq:rhat-assumptions}, and $F_{\rhat}(s,x)=0$, then
  $x\in sK(\eta,\rhat)$.  (The function $F_{\rhat}$ is defined
  in~\eqref{eq:def-F}.)
\end{lemma}
\begin{proof}
  Recall the functions $a$ and $b$ defined in
  Proposition~\ref{prop:F-properties}. Note that by
  inequalities~\eqref{eq:a-estimate} and~\eqref{eq:b-estimate}, for
  every $x\in\R^2\setminus\{0\}$ the inequality
  \begin{equation*}
    a(x)^2+b(x)^2 < 2\mu^2
  \end{equation*}
  holds, and that it is possible to find a constant $L_1=L_1(\eta)>0$
  so that $|x|\ge L_1$ implies
  \begin{equation}
    \label{eq:i-of-large-ineq}
    \frac{1}{4} \le \frac{1}{a(x)^2+b(x)^2} \le \frac{9}{4}
  \end{equation}
  and
  \begin{equation}
    \label{eq:i-of-large-eta}
    \frac{a(x)^2}{a(x)^2+b(x)^2} \ge \eta^2.
  \end{equation}

  Now if $x\in\R^2\setminus\{0\}$ satisfies $F_{\rhat}(s,x)=0$, that
  is, $X(y(x))x=s\rhat$, then
  \begin{equation}
    \label{eq:i-of-large-Pythagoras}
    s^2 = |x|^2(a(x)^2+b(x)^2) < 2\mu^2|x|^2.
  \end{equation}
  Therefore if $s\ge\sqrt{2}\mu L_1$ and $F_{\rhat}(s,x)=0$, then
  $|x|>L_1$, and by~\eqref{eq:i-of-large-ineq}
  and~\eqref{eq:i-of-large-Pythagoras}
  \begin{equation*}
    \frac{1}{2} \le \frac{|x|}{s} \le \frac{3}{2},
  \end{equation*}
  and by~\eqref{eq:i-of-large-ineq} and~\eqref{eq:i-of-large-eta}
  \begin{equation*}
    x\cdot\rhat
    = \frac{|x|}{s}\,\xhat\cdot X(y(x))x
    = \frac{|x|a(x)}{\sqrt{a(x)^2+b(x)^2}}
    \ge \eta|x|.
  \end{equation*}
  It follows that $x/s\in K(\eta,\rhat)$, and consequently the lemma
  holds if $L\ge \sqrt{2}\mu L_1$.
\end{proof}

Given $s>0$ and $\rhat\in\R^2$, define a mapping
$G_{s\rhat}:\R^2\to\R^2$ by
\begin{equation}
  \label{eq:G-def}
  G_{s\rhat}(x) := x-F_{\rhat}(s,x) = x - X(y(x))x + s\rhat.
\end{equation}
Obviously, for every $s>0$ the set of zeros of $F_{\rhat}(s,\cdot)$
and the set of fixed points of $G_{s\rhat}$ coincide.

\begin{lemma}
  Let $0<\eta<1$. There exists a constant $L=L(\eta)>0$ such that if
  $s\ge L$ and $\rhat\in\R^2$ is a vector that
  satisfies~\eqref{eq:rhat-assumptions}, then
  \begin{equation}
    \label{eq:K-inclusion}
    G_{s\rhat}\big[sK(\eta,\rhat)\big]\subset sK(\eta,\rhat).
  \end{equation}
\end{lemma}
\begin{proof}
  Let ${w,\rhat}\in\R^2$ satisfy $1/2\le|w|\le3/2$ and $|\rhat|=1$.
  With the notation of Proposition~\ref{prop:F-properties}, for $s>0$,
  \begin{equation*}
    \frac{1}{s}G_{s\rhat}(sw) - \rhat
    = \big(1-a(sw)\big)w - b(sw)|w|\widehat{w}_\perp
    := e(s,w,\rhat).
  \end{equation*}
  From inequalities~\eqref{eq:a-estimate} and~\eqref{eq:b-estimate} it
  follows $e(s,w,\rhat)\to0$ as $s\to\infty$, uniformly in $w$ and
  $\rhat$. It follows that there exists $L>0$ such that if $s\ge L$
  and $w\in K(\eta,\rhat)$, then $G_{s\rhat}(sw)/s\in
  K(\eta,\rhat)$. This implies~\eqref{eq:K-inclusion}.
\end{proof}

Below $D_xG_{s\rhat}$ denotes the Jacobian matrix of the map
$G_{s\rhat}$ defined in~\eqref{eq:G-def}.
\begin{lemma}
  \label{lemma:small-D}
  Let $\rhat\in\R^2$ satisfy~\eqref{eq:rhat-assumptions}. There exists
  numbers $\eta = \eta(\rhat)\in(0,1)$ and $L=L(\rhat)>0$ such that if
  $s\ge L$, ${x,x'}\in sK(\eta,\rhat)$, and $0\le\nu\le 1$, then
  \begin{equation}
    \label{eq:small-D}
    \big\|D_xG_{s\rhat}((1-\nu)x+\nu x')\big\| \le \frac{1}{2}.
  \end{equation}
  Here the norm is the operator norm on $\R^{2\times 2}$.
\end{lemma}
\begin{proof}
  An expression for $D_xG_{s\rhat}(x)$ is readily obtained from that
  of $D_xF_{\rhat}(x)$, which was calculated
  in~\eqref{eq:DxF-expression}. Observe that all of the polynomials
  $p_{ij}$ in~\eqref{eq:DxF-expression} have total degrees at most
  six.

  Let $C>0$ be large enough so that
  $\|D_xG_{s\rhat}(x)\|\le C|x|^6/\det Y(x)^2$ for every $x\in\R^2$
  with $|x|\ge 1$. Next, choose a constant $\eta=\eta(\rhat)\in(0,1)$
  so that if $x\in\R^2$ and $x\cdot\rhat\ge\eta|x|$, then
  $x_1\ge|x|\rhat_1/\sqrt{2}$ and $x_2\ge|x|\rhat_2/\sqrt{2}$. With
  these constants, for every $x\in\R^2$ with $x\cdot\rhat\ge\eta|x|$
  and $|x|\ge 1$, it holds that
  \begin{equation}
    \label{eq:small-D-ineq}
    \big\|DG_{s\rhat}(x)\big\|
    \le \frac{C|x|^6}{(\rhat_1\rhat_2)^4\, |x|^8}=\frac{C'}{|x|^2},
  \end{equation}
  where $C':=C/(\rhat_1\rhat_2)^4>0$.

  Now consider $x=sw$ and $x'=sw'$, where $s>0$ and
  ${w,w'}\in K(\eta,\rhat)$. If $0\le\nu\le 1$, then for
  $w_\nu:=(1-\nu)w+\nu w'$ both of the inequalities $|w_\nu|^2\ge 1/8$
  and $w_\nu\cdot\rhat\ge\eta|w_\nu|$ hold. Therefore, if $s^2\ge 8$,
  it follows from~\eqref{eq:small-D-ineq} that
  \begin{equation*}
    \big\|D_xG_{s\rhat}((1-\nu)x+\nu x')\big\|
    = \big\|D_xG_{s\rhat}(sw_\nu)\big\|
    \le \frac{C'}{|sw_\nu|^2}
    \le \frac{8C'}{s^2}.
  \end{equation*}
  Consequently, for $s\ge\max\{2\sqrt{2},4\sqrt{C'}\}$
  inequality~\eqref{eq:small-D} holds.
\end{proof}

\begin{proposition}
  \label{prop:fixed-point}
  Let $\rhat\in\R^2$ satisfy~\eqref{eq:rhat-assumptions}. There exists
  a constant $L=L(\rhat)>0$ such that following hold:
  \begin{enumerate}[(i)]
  \item For every $s\ge L$ the function $G_{s\rhat}$ has a unique
    fixed point in $\R^2$.
  \item If $h:[L,\infty)\to\R^2$ denotes the function that maps $s$ to
    the unique fixed point of $G_{s\rhat}$, then $h$ is differentiable
    on $(L,\infty)$.
  \item There exists a constant $C>0$ (independent of $\rhat$) such
    that the function $h$ from {(ii)} satisfies
    \begin{equation}
      \label{eq:xs-estimate}
      h(s) = s(\rhat+e(s)),\text{ where } |e(s)|\le\frac{C}{s^{2/3}}.
    \end{equation}
  \end{enumerate}
\end{proposition}

\begin{proof}
  Let $\eta=\eta(\rhat)\in(0,1)$ and $L=L(\rhat)>0$ be such that for
  $s\ge L$ inequality~\eqref{eq:small-D} holds for every
  ${x,x'}\in sK(\eta,\rhat)$ and $0\le\nu\le 1$. If necessary,
  increase $L$ so that in addition for $s\ge L$
  inclusion~\eqref{eq:K-inclusion} holds and equality
  $F_{\rhat}(s,x)=0$ implies that $x\in sK(\eta,\rhat)$ (cf.\
  Lemma~\ref{lemma:i-of-large-lemma}).

  Let $s\ge L$. Then $G_{s\rhat}$ maps $sK(\eta,\rhat)$ into itself,
  and if ${x,x'}\in sK(\eta,\rhat)$, applying the fundamental theorem
  of calculus and estimating with~\eqref{eq:small-D} shows that
  \begin{equation*}
    |G_{s\rhat}(x)-G_{s\rhat}(x')|
    \le |x-x'|\sup_{0\le\nu\le 1}\big\|D_xG_{s\rhat}((1-\nu)x+\nu x')\big\|
    \le\frac{|x-x'|}{2}.
  \end{equation*}
  Thus, the restriction of $G_{s\rhat}$ to $sK(\eta,\rhat)$ is a
  contraction.

  By the Banach fixed-point theorem the function $G_{s\rhat}$ has a
  unique fixed point in $sK(\eta,\rhat)$. Because $G_{s\rhat}(x)=x$ if
  and only if $F_{\rhat}(s,x)=0$, this fixed point is unique in
  $\R^2$, also. Part {(i)} is now proved.
    
  Let $s_0>L$ and $h$ be as in {(ii)}. Consider the function
  $(L,\infty)\times\R^2\ni(s,x)\mapsto F_{\rhat}(s,x)\in\R^2$ at a
  neighborhood of its zero $(s_0,h(s_0))$. Since
  \begin{equation*}
    D_xF_{\rhat}(x) = I_2 - D_xG_{s\rhat}(x),
  \end{equation*}
  it follows from inequality~\eqref{eq:small-D} that at the point
  $(s,x)=(s_0,h(s_0))$ the derivative $D_xF_{\rhat}(x)$ is
  invertible. Then by the implicit function theorem in some
  neighborhood $(s_0-\epsilon,s_0+\epsilon)$ the zero of
  $F_{\rhat}(s,\cdot)$, i.e., $h(s)$, depends differentiably on
  $s$. Because $s_0>L$ was arbitrary, the function $s\mapsto h(s)$ is
  differentiable, and {(ii)} is proved.

  If we write $h(s)$ as in~\eqref{eq:xs-estimate} and denote
  $x:=h(s)$, then in the notation of
  Proposition~\ref{prop:F-properties} we have
  \begin{equation*}
    e(s)
    = \frac{1}{s}x-\rhat
    = \frac{1}{s}G_{s\rhat}(x) - \rhat
    = \frac{|x|}{s}\big[\big(1-a(x)\big)\xhat-b(x)\xhat_{\perp}\big].
  \end{equation*}
  Because $1/2\le|x|/s\le3/2$ since $x\in sK(\eta,\rhat)$, we obtain
  from~\eqref{eq:a-estimate} and~\eqref{eq:b-estimate} that for some
  constant $C>0$ depending only on $\mu$ it holds that
  $|e(s)|\le C/s^{2/3}$, for every $s\ge L$. This proves {(iii)}.
\end{proof}

With the previous proposition in hand, we can now prove the main
theorem of this section. Note that, among others, the theorem states
that unlike in the case of weak injected fields, in which case
system~\eqref{eq:system} has nine equilibrium points
(Theorem~\ref{thm:equilibrium-small-dynamics}), in the case of strong
injected fields, the system has a single equilibrium point.

\begin{theorem}
  \label{thm:strong-injection}
  Consider $\uhat=(\uhat_-, \uhat_+)\in\C^2$ with $\uhat_-\neq 0$ and
  $\uhat_+\neq 0$. There exists a constant $L=L(\uhat)>0$ and a
  continuous function
  \begin{equation*}
    E_{\uhat} : \{\lambda\in\C : |\lambda|\ge L\}\to\C^2
  \end{equation*}
  with the following property: If in system~\eqref{eq:system} the
  injected field $u$ is of the form $u=\lambda\uhat$ with
  $|\lambda|\ge L$, then a triple $(E, N, n)\in\C^2\times\R\times\R$
  is an equilibrium point of the system, if and only if
  \begin{equation*}
    E = E_{\uhat}(\lambda) \text{ and } (N,n) = y(|E_-|,|E_+|) 
  \end{equation*}
  (the function $y$ is defined in~\eqref{eq:def-y}). Furthermore,
  there exists a constant $C=C(\uhat)>0$ such that the function
  $E_{\uhat}$ satisfies
  \begin{equation}
    \label{eq:Eu-estimate}
    E_{\uhat}(\lambda)
    = \frac{\lambda e^{i\theta}}{|1+i\alpha|}(\uhat+e(\lambda))
    \text{, where } |e(\lambda)| \le\frac{C}{|\lambda|^{2/3}}
    \text{ and }\theta:=-\arg(1+i\alpha).
  \end{equation}
\end{theorem}

\begin{remark}
  It follows from~\eqref{eq:Eu-estimate} that the magnitudes of the
  emitted field $E_{\uhat}(\lambda)$ and the injected field
  $u=\lambda\uhat$ are asymptotically related by
  \begin{equation*}
    \lim_{|\lambda|\to\infty}\frac{|E_{\uhat}(\lambda)|}{|\lambda\uhat|}
    = \frac{1}{|1+i\alpha|},
  \end{equation*}
  and that as $\lambda$ grows, the polarization of the emitted field
  $E_{\uhat}(\lambda)$ approaches on the normalized Poincar\'e sphere
  that of $\uhat$.
\end{remark}

\begin{proof}
  Define
  \begin{equation}
    \label{eq:rhat-strong}
    \rhat := \frac{1}{|\uhat|}
    \begin{bmatrix}
      |\uhat_-|\\|\uhat_+|
    \end{bmatrix}.
  \end{equation}
  Then $\rhat$ satisfies~\eqref{eq:rhat-assumptions}, let
  $L'=L'(\rhat)>0$ be a constant and $h:[L',\infty)\to\R^2$ a function
  as in Proposition~\ref{prop:fixed-point}.

  Fix a constant $L>|1+i\alpha||\uhat|\inv L'$, and define for
  $\lambda\in\C$ with $|\lambda|\ge L$ a function $E_{\uhat}$ by
  \begin{equation*}
    E_{\uhat}(\lambda)
    := e^{i\theta}\frac{\lambda}{|\lambda|}
    \begin{bmatrix}
      \frac{\uhat_-}{|\uhat_-|} & 0\\
      0 & \frac{\uhat_+}{|\uhat_+|}
    \end{bmatrix}
    h\left(\tfrac{|\lambda\uhat|}{|1+i\alpha|}\right).
  \end{equation*}
  As $h$ is differentiable on $(L',\infty)$, the function
  $E_{\uhat}(\lambda)$ is continuous on its domain. Also,
  estimate~\eqref{eq:Eu-estimate} follows directly
  from~\eqref{eq:xs-estimate}.

  Now with $s:=|1+i\alpha|\inv|\lambda\uhat|$ and $x:=h(s)$ it holds
  that $X(y(x))x=s\rhat$, so by
  Proposition~\ref{prop:algebraic-version} the triple
  $(E,N,n)\in\C^2\times\R\times\R$ with $E=E_{\uhat}(\lambda)$ and
  $(N,n)=y(x)=y(|x_1|,|x_2|)=y(|E_-|,|E_+|)$ is an equilibrium point
  of system~\eqref{eq:system} with injected field
  \begin{equation*}
    u = (1+i\alpha)
    e^{i\theta}\frac{\lambda}{|\lambda|}
    \begin{bmatrix}
      \frac{\uhat_-}{|\uhat_-|} & 0\\
      0 & \frac{\uhat_+}{|\uhat_+|}
    \end{bmatrix}
    s\rhat
    = \lambda\uhat.
  \end{equation*}

  On the other hand, consider an arbitrary equilibrium point $(E,N,n)$
  of system~\eqref{eq:system} with $u=\lambda\uhat$, where
  $|\lambda|\ge L$. By Proposition~\ref{prop:algebraic-version} there
  exists $x\in\R^2$, $s\ge 0$, $\rhat\in[0,\infty)\times[0,\infty)$
  with $|\rhat|$=1, and $\phi_\pm\in\R$ such that
  \begin{subequations}
    \begin{align}
      X(y(x))x &= s\rhat,\label{eq:necessary2-1}\\
      E &=
          \begin{bmatrix}
            x_1\,e^{i\phi_-}\\
            x_2\,e^{i\phi_+}
          \end{bmatrix}
      ,\label{eq:necessary2-2}\\
      \begin{bmatrix}
        N\\n
      \end{bmatrix}
               & = y(x),\text{ and}\label{eq:necessary2-3}\\
      \lambda\uhat &= (1+i\alpha)
                     \begin{bmatrix}
                       s\rhat_1\, e^{i\phi_-}\\
                       s\rhat_2\, e^{i\phi_+}
                     \end{bmatrix}
      \label{eq:necessary2-4}.
    \end{align}
  \end{subequations}
  
  Equation~\eqref{eq:necessary2-4} implies that
  $s=|1+i\alpha|\inv|\lambda\uhat|>L'$ and that $\rhat$
  satisfies~\eqref{eq:rhat-strong}. Then from~\eqref{eq:necessary2-1}
  it follows that $G_{s\rhat}(x)=x$, so $x=h(s)$ by
  Proposition~\ref{prop:fixed-point}. The numbers $e^{i\phi_\pm}$ can
  be determined from~\eqref{eq:necessary2-4}, inserting them
  into~\eqref{eq:necessary2-2} shows that
  $E=E_{\uhat}(\lambda)$. Finally, from~\eqref{eq:necessary2-3}
  and~\eqref{eq:necessary2-2} it follows that $(N,n)=y(|E_-|,|E_+|)$.
\end{proof}

%
\section{Optical neural networks based on injection locking}
\label{sec:neural-network}

\begin{figure}[t]
  \hspace{-1.25em}
  \begin{subfigure}[b]{0.54\textwidth}
    \centering \scalebox{0.9}{\begin{tikzpicture}


    \draw[draw=black, line width=1.2pt, fill=cyan!30, rounded corners=5mm] (-0.5,-0.5) rectangle (0.5,5.5);
    \draw[draw=black, line width=1.2pt, fill=purple!30, rounded corners=5mm] (2.5,-0.5) rectangle (3.5,5.5);
    \draw[draw=black, line width=1.2pt, fill=cyan!30, rounded corners=5mm] (5.5,-0.5) rectangle (6.5,5.5);
    \node[above, align=center] at (0,5.7) {Input layer\\$I$};
    \node[above, align=center] at (3,5.7) {Hidden layer\\$J$};
    \node[above, align=center] at (6,5.7) {Output layer\\$K$};

    \begin{scope}[thick,decoration={
        markings,
        mark=at position 0.37 with {\arrow{Latex[round]}}}
        ] 
        \draw[postaction={decorate}] (0,0) -- (3, 0);
        \draw[postaction={decorate}] (0,0) -- (3, 2);
        \draw[postaction={decorate}] (0,2) -- (3, 0);
        \draw[postaction={decorate}] (0,2) -- (3, 2);
        \draw[postaction={decorate}] (0,2) -- (3, 4);
        \draw[postaction={decorate}] (0,4) -- (3, 2);
        \draw[postaction={decorate}] (0,4) -- (3, 4);
        \draw[postaction={decorate}] (0,4) -- (3, 5);
        \draw[postaction={decorate}] (0,5) -- (3, 4);
        \draw[postaction={decorate}] (0,5) -- (3, 5);
        \draw[postaction={decorate}] (3,0) -- (6, 0);
        \draw[postaction={decorate}] (3,0) -- (6, 2);
        \draw[postaction={decorate}] (3,2) -- (6, 0);
        \draw[postaction={decorate}] (3,2) -- (6, 2);
        \draw[postaction={decorate}] (3,2) -- (6, 4);
        \draw[postaction={decorate}] (3,4) -- (6, 2);
        \draw[postaction={decorate}] (3,4) -- (6, 4);
        \draw[postaction={decorate}] (3,4) -- (6, 5);
        \draw[postaction={decorate}] (3,5) -- (6, 4);
        \draw[postaction={decorate}] (3,5) -- (6, 5);
    \end{scope}
    \begin{scope}[line width=1.8pt, blue!60, decoration={
        markings,
        mark=at position 0.75 with {\arrow{Latex}}}
        ] 
        \draw[postaction={decorate}] (2.3,0.7) -- (3, 0);
        \draw[postaction={decorate}] (2.3,2.7) -- (3, 2);
        \draw[postaction={decorate}] (2.3,4.7) -- (3, 4);
        \draw[postaction={decorate}] (2.3,5.7) -- (3, 5);
    \end{scope}
    \begin{scope}[thick, decoration={
        markings,
        mark=at position 1.0 with {\arrow{Latex[round]}}}
        ] 
        \draw[postaction={decorate}] (6,0) -- (7,0);
        \draw[postaction={decorate}] (6,2) -- (7,2);
        \draw[postaction={decorate}] (6,4) -- (7,4);
        \draw[postaction={decorate}] (6,5) -- (7,5);
    \end{scope}

    \draw[draw=gray, line width=1.2pt, fill=cyan!10] (0,5) circle [radius=0.35] node {$1$};
    \draw[draw=gray, line width=1.2pt, fill=cyan!10] (0,4) circle [radius=0.35] node {$2$};
    \draw[draw=none] (0,3) circle [radius=0.35] node {$\vdots$};
    \draw[draw=gray, line width=1.2pt, fill=cyan!10] (0,2) circle [radius=0.35] node {$N_i$};
    \draw[draw=none] (0,1) circle [radius=0.35] node {$\vdots$};
    \draw[draw=gray, line width=1.2pt, fill=cyan!10] (0,0) circle [radius=0.35] node {$i_0$};

    \draw[draw=gray, line width=1.2pt, fill=purple!10] (3,5) circle [radius=0.35] node {$1$};
    \draw[draw=gray, line width=1.2pt, fill=purple!10] (3,4) circle [radius=0.35] node {$2$};
    \draw[draw=none] (3,3) circle [radius=0.35] node {$\vdots$};
    \draw[draw=gray, line width=1.2pt, fill=purple!10] (3,2) circle [radius=0.35] node {$N_j$};
    \draw[draw=none] (3,1) circle [radius=0.35] node {$\vdots$};
    \draw[draw=gray, line width=1.2pt, fill=purple!10] (3,0) circle [radius=0.35] node {$j_0$};

    \draw[draw=gray, line width=1.2pt, fill=cyan!10] (6,5) circle [radius=0.35] node {$1$};
    \draw[draw=gray, line width=1.2pt, fill=cyan!10] (6,4) circle [radius=0.35] node {$2$};
    \draw[draw=none] (6,3) circle [radius=0.35] node {$\vdots$};
    \draw[draw=gray, line width=1.2pt, fill=cyan!10] (6,2) circle [radius=0.35] node {$N_k$};
    \draw[draw=none] (6,1) circle [radius=0.35] node {$\vdots$};
    \draw[draw=gray, line width=1.2pt, fill=cyan!10] (6,0) circle [radius=0.35] node {$k_0$};

    \node[above, align=left] at (1.5,2) {$E_i^{\left( I \right) }$};
    \node[above, align=left] at (4.5,2) {$E_j^{\left( J \right) }$};
    \node[above, align=left] at (7.22,2) {$E_k^{\left( K \right) }$};
    \node[above, align=right] at (2.4,2.7) {$E_j^{\left( \textrm{ext} \right) }$};

\end{tikzpicture}}
    \caption{Optical neural network}
    \label{fig:neural-network}
  \end{subfigure}
  \hspace{1.5em}
  \begin{subfigure}[b]{0.45\textwidth}
    \centering \scalebox{0.95}{\input{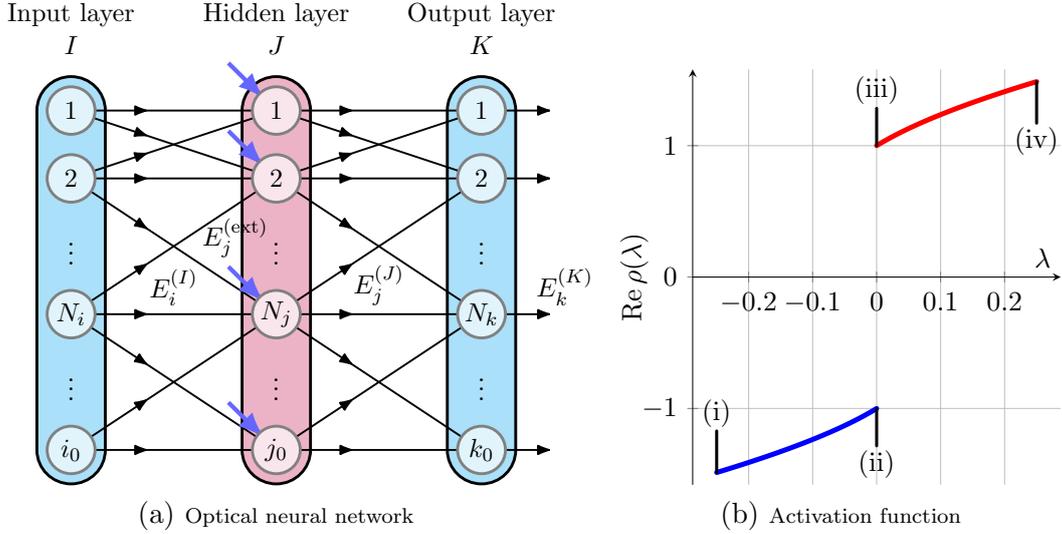}}
    \caption{Activation function}
    \label{fig:activation-function}
  \end{subfigure}
  \caption[Optical neural network]{Schematic illustration of an
    optical neural network and a complex-valued activation function
    $\rho$ based on injection locking. The parameters in {(b)} are
    those of Figures~\ref{fig:ODE-solution} to~\ref{fig:stability}. In
    this figure, polarization $\uhat$ of the electric fields in the
    network has been chosen so that $\im\rho(\lambda)=0$ for
    $\lambda\in\R$. Labels {(i)}--{(iv)} in {(b)} match those of
    Figures~\ref{fig:paths} and~\ref{fig:stability}.
    
    Fields $E_i^{(I)}=\lambda_i^{(I)}\uhat$ in the input layer $I$ are
    inputs to the network. They are passed through passive optical
    elements (which correspond to multiplication by $a_{ji}\in\C$) and
    joined with fixed external fields $E_j^{(\mathrm{ext})}=b_j\uhat$
    to form a field
    $(\sum_i a_{ji}\lambda_i^{(I)}+b_j)\uhat = \lambda_j^{(J)}\uhat$
    injected into the $j$:th laser in the hidden layer $J$. Due to
    injection locking, the corresponding emitted field $E_j^{(J)}$ is
    $\rho(\lambda_j^{(J)})\uhat$. The fields from the hidden layer are
    passed through passive optical elements and joined to form outputs
    $E_k^{(K)} = \lambda_k^{(K)}\uhat$ of the network.}
  \label{fig:neural-network-ab}
\end{figure}

We now describe a design of an optical neural network that can be
implemented with a network of lasers, and whose working principle is
based on injection locking (see Figure~\ref{fig:neural-network}). The
network consists of an input layer (Layer~$I$), an output layer
(Layer~$K$), and one hidden layer (Layer~$J$) in between (the working
principle naturally generalizes to a network with several hidden
layers):
\begin{enumerate}[(i)]
\item In the input layer, each node (artificial neuron) is a
  laser. The nodes in this layer are not connected to each other, and
  the output of a node is the electric field emitted by the
  corresponding laser.
\item In the hidden layer, the nodes are lasers that are coupled to
  injected electric fields. The injected fields are composed of fixed
  external electric fields together with outputs of the input layer
  modified by some passive optical elements, e.g., polarizers or
  mirrors, optical isolators, and absorbing components. Due to
  injection locking, each laser in the hidden layer stabilizes to some
  equilibrium point determined by the injected field, and the output
  of a node is the emitted electric field.

  The coupling between layers $I$ and $J$ is unidirectional, we note
  that one can use lasers of varying powers to replace the use of
  optical isolators.
\item Between the hidden layer and the output layer, the electric
  fields from the hidden layer are first modified by passive optical
  elements, and then joined to form the output of the network. The
  nodes in the output layer correspond to exits of optical cables or
  waveguides in integrated optics.
\end{enumerate}.

The relation between inputs and outputs of the network is set by
choosing the external electric fields that are part of the injected
fields in the hidden layer, and the passive optical elements on both
sides of the hidden layer. We will show that an arbitrary continuous
function can be approximated within any given accuracy by networks of
this form.

The optical neural network is modeled mathematically as follows.
Indexes of lasers in the input layer are denoted by
$I=\{1,2,\dots,i_0\}$. The output of $i$:th laser is a linearly
polarized electric field $E_i^{(I)}\in\C^2$, and all electric fields
in this layer are assumed to share the same linear polarization, i.e.,
for all $i=1,2,\ldots,i_0$,
\begin{equation}
  \label{eq:EJ}
  E_i^{(I)} = \lambda_i^{(I)}\uhat,
\end{equation}
where $\lambda_i^{(I)}\in\C$, and
$\uhat=(\uhat_-, \uhat_+)\in\C^2\setminus\{0\}$ is fixed and satisfies
$|\uhat_-|=|\uhat_+|$. It is also assumed that the set of all possible
inputs is bounded, i.e., there exists $R>0$ such that whenever
$(\lambda_i^{(I)}\uhat)_{i=1}^{i_0}$ is an input to the network, then
$|(\lambda_i^{(I)})_{i=1}^{i_0}|_{\C^{i_0}}\le R$. Here
$|\cdot|_{\C^{i_0}}$ denotes the Euclidean norm on $\C^{i_0}$.

In the hidden layer indexes of lasers are denoted by
$J=\{1,2,\dots,j_0\}$. The passive optical elements between the input
layer and the hidden layer may induce scaling and phase shift to the
electric fields, i.e., field $E_i^{(I)}$ from the $i$:th laser of the
input layer to the $j$:th laser of the hidden layer transforms to
$a_{ji}E_i^{(I)}$, where $a_{ji}\in\C$. The total injected field
$u_j\in\C^2$ to the $j$:th laser in the hidden layer is then the sum
of the modified fields and an external electric field
$E_j^{\mathrm{(ext)}}$, which is assumed to share the same
polarization with the lasers in the input layer:
$E_j^{\mathrm{(ext)}}=b_j\uhat$ for some $b_j\in\C$. Thus,
\begin{equation}
  u_j
  = \sum_{i=1}^{i_0} a_{ji} E_i^{(I)} + E_j^{\mathrm{(ext)}}
  = \left(\sum_{i=1}^{i_0} a_{ji}\lambda_i^{(I)} + b_j\right)\uhat.
\end{equation}

By Theorems~\ref{thm:equilibrium-small-dynamics}
and~\ref{thm:stability}, if the linewidth enhancement factor $\alpha$
of the laser is zero (i.e., $\alpha=0$ in system~\eqref{eq:system})
and the injected field $u_j$ to the $j$:th laser is written as
$u_j=\lambda_j^{(J)}\uhat$, then for some constant $\ell>0$ it holds
that as long as $0<|\lambda_j^{(J)}|<\ell$, then the $j$:th laser has
a unique stable equilibrium point (denoted by
$E_{\uhat}^{(+\textsc{x})}(\lambda_j^{(J)})$ in
Theorems~\ref{thm:equilibrium-small-dynamics}
and~\ref{thm:stability}). If $\alpha>0$, then this point is still an
equilibrium point, and it was shown in Section~\ref{sec:weak-fields}
how to numerically check if for weak enough injected fields it is a
unique stable equilibrium point.  Assuming this is the case, after a
successful injection locking the emitted field $E_j^{(J)}\in\C^2$ of
the $j$:th laser in the hidden layer with small enough injected field
$u_j=\lambda_j^{(J)}\uhat\neq 0$ stabilizes to
\begin{equation*}
  E_j^{(J)}
  = \rho(\lambda_j^{(J)})\uhat,  
\end{equation*}
where the function
\begin{equation}
  \label{eq:actfn}
  \rho:=\rho^{(+\textsc{x})}:\{\lambda\in\C : 0<|\lambda|<\ell\} \to \C
\end{equation}
is defined in Theorem~\ref{thm:equilibrium-small-dynamics}.
Figure~\ref{fig:neural-network-ab}b illustrates the function $\rho$
corresponding to the system in Figure~\ref{fig:ODE-solution}.

In the output layer nodes are indexed by $K=\{1,2,\ldots,k_0\}$, and
the $k$:th output $E_k^{(K)}\in\C^2$ of the network is a superposition
of the emitted fields $E_j^{(J)}$ of lasers in the hidden layer
modified by passive optical elements represented by complex numbers
$c_{kj}$:
\begin{equation}
  \label{eq:El}
  E_k^{(K)}
  = \sum_{j=1}^{j_0} c_{kj} E_j^{(J)}
  = \sum_{j=1}^{j_0} c_{kj}
  \rho(\lambda_j^{(J)})\uhat,
\end{equation}
whenever $0<|\lambda_j^{(J)}|<\ell$ for all $j=1,2,\ldots,j_0$.

As the input to the network is of the form
$(\lambda_i^{(I)}\uhat)_{i=1}^{i_0}\in(\C^2)^{i_0}$,
$\lambda_i^{(I)}\in\C$, and the output is by~\eqref{eq:El} of the form
$(\lambda_k^{(K)}\uhat)_{k=1}^{k_0}\in(\C^2)^{k_0}$,
$\lambda_k^{(K)}\in\C$, the network essentially computes the map
\begin{equation*}
  (\lambda^{(I)}_i)_{i=1}^{i_0}
  \mapsto
  (\lambda^{(K)}_k)_{k=1}^{k_0}
  =:\M((\lambda_i^{(I)})_{i=1}^{i_0}).
\end{equation*}
It follows from equations~\eqref{eq:EJ}--\eqref{eq:El} that the $k$:th
component function $\M_k$ of $\M$ is
\begin{equation}
  \label{eq:M-l}
  \M_k((\lambda_i^{(I)})_{i=1}^{i_0})
  = \sum_{j=1}^{j_0}c_{kj}\rho\left(\sum_{i=1}^{i_0} a_{ji}\lambda^{(I)}_i + b_j\right),
\end{equation}
where it is assumed that
\begin{equation}
  \label{eq:ONN-assumption}
  0<\left|\sum_{i=1}^{i_0} a_{ji}\lambda^{(I)}_i + b_j\right|
  <\ell
  \text{ for every } j=1,2,\ldots,j_0.
\end{equation}

In~\eqref{eq:M-l} and~\eqref{eq:ONN-assumption} parameters
${a_{ji},c_{kj}}\in\C$ correspond to the passive optical elements
between the layers, and parameters $b_j\in\C$ correspond to the fixed
external electric fields.

\begin{remark}
  The lasers in the input layer are not connected with each other,
  yet, the formulation assumes that the phase differences remain
  constant at the equilibrium point. As known, all oscillatory signal
  sources, lasers included, fluctuate in phase. This drift will
  inevitably invalidate the assumption of the constant phase
  difference between two lasers unless they share a common reference
  (seed) signal. Therefore, a practical implementation of a
  laser-based optical neural network will require a common
  narrow-linewidth reference signal that is used to lock enough lasers
  in the network. At the bare minimum, all lasers of the first layer
  must be injected from the same source. The phase of the injected
  reference light may be controlled individually for each network
  node, but the natural fluctuations of the reference must be
  experienced equally among the injected lasers. This arrangement is
  not unlike the clock signal of a digital computer that is used to
  synchronize operations between individual circuits.
\end{remark}

Below $\Bbar\subset\C^{i_0}$ is the closed ball of radius $R$ centered
at the origin.
\begin{theorem}
  \label{thm:ONN}
  Fix integers $i_0>0$ and $k_0>0$ and a number $R>0$, let $\rho$ be
  as in~\eqref{eq:actfn}, and consider an arbitrary continuous
  function $f:\Bbar\to\C^{k_0}$. Let $\epsilon>0$. There exists an
  integer $j_0>0$ and numbers ${a_{ji},b_j,c_{kj}}\in\C$,
  $j=1,2,\ldots,j_0$, $i=1,2,\ldots,i_0$, $k=1,2,\ldots,k_0$, such
  that following holds:
  \begin{enumerate}[(i)]
  \item The inequalities~\eqref{eq:ONN-assumption} hold for a.e.\
    $(\lambda_i^{(I)})_{i=1}^{i_0}\in\Bbar$ (the measure on
    $\Bbar\subset\C^{i_0}=\R^{2i_0}$ is the $2i_0$-dimensional
    Lebesgue measure), and
  \item the function $\M$ defined componentwise a.e.\ in $\Bbar$
    by~\eqref{eq:M-l} is measurable and satisfies
    \begin{equation}
      \label{eq:ONN}
      \big\| \M - f \big\|_{L^\infty(\Bbar;\C^{k_0})} \le \epsilon.
    \end{equation}
  \end{enumerate}
\end{theorem}
\begin{proof}
  Let $U:=\{\lambda\in\C:|\lambda|<\ell\}$ and extend the function
  $\rho$ defined in~\eqref{eq:actfn} into a function $\rho:U\to\C$ by
  setting $\rho(0):=0$. Then $\rho$ is locally bounded on $U$ and
  continuous on $U\setminus\{0\}$, and by
  Theorem~\ref{thm:equilibrium-small-dynamics}
  \begin{equation*}
    \lim_{\substack{\lambda\in\R,\\\lambda\to 0^+}}\rho(\lambda)
    = -\lim_{\substack{\lambda\in\R,\\\lambda\to 0^-}}\rho(\lambda)
    \neq 0.
  \end{equation*}
  In particular $\rho$ is not a.e.\ equal to a continuous function,
  and consequently it satisfies both {(i)} and {(ii)} of
  Theorem~\ref{thm:UAT} stated in the Appendix (note that if
  $\Delta^m\rho\equiv 0$ for some $m\in\N$ in the sense of
  distributions, then $\rho$ is a.e.\ equal to a smooth function by
  elliptic regularity~\cite{MR1157815}).

  Let $f:\Bbar\to\C^{k_0}$ be a continuous function and fix
  $\epsilon>0$. By Theorem~\ref{thm:UAT} there exists an integer
  $j_0>0$ and parameters ${a_{ji},b_j,c_{kj}}\in\C$ such that
  \begin{equation*}
    \sum_{i=1}^{i_0} a_{ji}\lambda_i+b_j \in U
  \end{equation*}
  for every $j=1,2,\ldots,j_0$ and $(\lambda_i)_{i=1}^{i_0}\in\Bbar$,
  and such that the network $\mathcal{N}:\Bbar\to\C^{k_0}$ defined
  componentwise by~\eqref{eq:C-NN} satisfies
  \begin{equation*}
    \sup_{(\lambda_i)\in\Bbar} \big| \mathcal{N}((\lambda_i)_{i=1}^{i_0}) - f((\lambda_i)_{i=1}^{i_0}) \big|_{\C_{k_0}} \le\epsilon.
  \end{equation*}
  Furthermore, it may be assumed that for every $j$ either
  $(a_{j1},a_{j2},\ldots,a_{ji_0})\neq 0$ or $b_j\neq 0$, since
  otherwise the corresponding term does not affect the value of
  $\mathcal{N}$. Observe that $\mathcal{N}$ is measurable, because the
  set
  \begin{equation*}
    N := \bigcup_{j=1}^{j_0}
    \Big\{(\lambda_i)_{i=1}^{i_0}\in\C^{i_0} :
    \sum_{i=1}^{i_0 }a_{ji}\lambda_i + b_j = 0\Big\}
  \end{equation*}
  has $2i_0$-dimensional Lebesgue measure zero and the restriction of
  $\mathcal{N}$ to $\Bbar\setminus N$ is continuous.

  Let us define $\M$ by the same parameters $j_0$, $a_{ji}$, $b_j$ and
  $c_{kj}$ as $\mathcal{N}$. Because
  inequalities~\eqref{eq:ONN-assumption} hold on $\Bbar\setminus N$,
  the function $\M$ is defined a.e.\ in $\Bbar$. Furthermore,
  $\M=\mathcal{N}$ a.e.\ in $\Bbar$, so $\M$ is measurable and
  inequality~\eqref{eq:ONN} holds.
\end{proof}

%

\section*{Acknowledgment}

ML and LY were supported by the Academy of Finland (Finnish Centre of
Excellence in Inverse Modelling and Imaging and projects 273979,
284715, and 312110).

\appendix
\section*{Appendix: Approximation theorem for complex-valued neural
  networks}

In this appendix, we generalize the recent universal approximation
theorem for complex-valued neural networks by
F.~Voigtlaender~\cite{voigtlaender2020universal} to the case of
activation functions defined locally in an open subset $U\subset\C$,
instead of globally on the whole complex plane. The gist of the proof,
namely the use of Wirtinger calculus~\cite{MR716497} to show that the
functions $z^\alpha\zbar^\beta$ ($\zbar$ is the complex conjugate of
$z$) can be approximated by neural networks, is the same as in the
proof of Voigtlaender's theorem. However, the proof is complicated by
the fact that parameters for the network need to be chosen so that all
inputs to the activation function stay within $U$.

Let $\Bbar:=\{z\in\C^{i_0}:|z|_{\C^{i_0}}\le R\}$. We consider
(shallow) complex-valued neural networks
$\mathcal{N}:\Bbar\to\C^{k_0}$, whose $k$:th component function is of
the form
\begin{equation}
  \label{eq:C-NN}
  \mathcal{N}_k(z) := \sum_{j=1}^{j_0}c_{kj}\rho(a_j\cdot z + b_j),
\end{equation}
where $a_j\cdot z := \sum_i a_{ji}z_i$. Here the integers $i_0>0$,
$j_0>0$, and $k_0>0$ are the number of inputs of the network, the
width of the network, and the number of outputs of the network,
respectively, and $\rho:U\to\C$, where $U\subset\C$ is an open set, is
the activation function. The parameters
$a_j=(a_{j1}, a_{j2},\ldots,a_{ji_0})\in\C^{i_0}$, $j=1,2,\ldots,j_0$,
$b\in\C^{j_0}$, and $(c_{kj})\in\C^{k_0\times j_0}$ are required to
satisfy
\begin{equation}
  \label{eq:parameter-requirement}
  a_j\cdot z + b_j\in U\text{ for every } z
  \in\Bbar\text{ and } j=1,2,\ldots,j_0.
\end{equation}

Following theorem is a local version of Voigtlaender's universal
approximation theorem for complex-valued neural
networks~\cite[Theorem~{1.3}]{voigtlaender2020universal}:
\begin{theorem}
  \label{thm:UAT}
  Let $i_0$, $k_0$, $R$, and $\rho$ be as above, and suppose that
  \begin{enumerate}[(i)]
  \item $\rho$ is locally bounded and continuous almost everywhere in
    the nonempty open set $U\subset\C=\R^2$ (the measure is the
    two-dimensional Lebesgue measure), and
  \item $\Delta^m\rho$ does not vanish identically in $U$ for any
    $m=0,1,2,\ldots$ (here
    $\Delta=\partial^2/\partial x^2+\partial^2/\partial y^2$,
    $z=x+iy$, is the Laplace operator defined in the sense of
    distributions).
  \end{enumerate}
  If $f:\Bbar\to\C^{k_0}$ is continuous and $\epsilon>0$, then there
  exists an integer $j_0>0$ and parameters $a_j\in\C^{i_0}$,
  $j=1,2,\ldots,j_0$, $b\in\C^{j_0}$, and
  $(c_{kj})\in\C^{k_0\times j_0}$ such
  that~\eqref{eq:parameter-requirement} holds, and that the
  complex-valued neural network $\mathcal{N}$ defined componentwise
  by~\eqref{eq:C-NN} satisfies
  \begin{equation}
    \label{eq:universal-approximation}
    \sup_{z\in\Bbar} \big|
    \mathcal{N}(z) - f(z)
    \big|_{\C^{k_0}}
    \le \epsilon.
  \end{equation}
\end{theorem}

There is a slight difference in the continuity assumption for the
activation function $\rho$ between Theorem~\ref{thm:UAT}
and~\cite[Theorem~{1.3}]{voigtlaender2020universal}. Here we require
that $\rho$ is continuous almost everywhere, i.e., that the set
$D\subset\C$ of its discontinuities is a null
set. In~\cite{voigtlaender2020universal} it is required that also the
closure of $D$ is a null set. The difference is due to how the
(potentially nonsmooth) activation function is smoothly approximated;
our approximation method is contained in the following two lemmas. Our
approach is similar to~\cite[Lemma~4]{hornik1993some}, in which
real-valued activation functions are considered. Theorem~\ref{thm:UAT}
will be proved after the lemmas.

\begin{lemma}
  \label{lemma:Riemann-type}
  For $\eta>0$, let $\PP(\eta)$ denote the set of countable partitions
  of $\R^2$ into measurable subsets with diameter at most $\eta$, and
  let $\psi:\R^2\to\C$ be a bounded and almost everywhere continuous
  function with compact support. Then
  \begin{equation}
    \label{eq:Riemann-like}
    \lim_{\eta\to 0}\left(
      \sup\left\{\sum_{j=1}^\infty\lambda_2(C_j)\,\sup_{{y,y'}\in C_j}|\psi(y)-\psi(y')|
        : (C_j)_{j=1}^\infty\in\PP(\eta)\right\}
    \right)
    = 0,
  \end{equation}
  where $\lambda_2$ denotes the Lebesgue measure on $\R^2$.
\end{lemma}
\begin{proof}
  Choose a sequence of partitions
  $((C_j(k))_{j=1}^\infty)_{k=1}^\infty\in\PP(1/k)$, and define
  \begin{equation*}
    d_k(x) := \sum_{j=1}^\infty
    \sup_{y,y'\in C_j(k)}|\psi(y)-\psi(y')|\, 1_{C_j(k)}(x),
  \end{equation*}
  where $1_{C_j(k)}$ is the characteristic function of the set
  $C_j(k)$.

  The functions $d_k$ are measurable, uniformly bounded by
  $2\|\psi\|_\infty$, and they are all supported in a fixed compact
  set. If $x\in\R^2$ is a point of continuity of $\psi$, then
  $d_k(x)\to 0$. As a consequence, $d_k\to 0$ as $k\to\infty$ almost
  everywhere in $\R^2$, and by the Lebesgue's dominated convergence
  theorem
  \begin{equation}
    \label{eq:Riemann-like-2}
    0
    = \lim_{k\to\infty} \int_{\R^2} d_k(x)\,dx
    = \lim_{k\to\infty}\left(
      \sum_{j=1}^\infty\lambda_2(C_j(k))\sup_{{y,y'}\in C_j(k)}|\psi(y)-\psi(y')|
    \right).
  \end{equation}
  This proves the lemma as the sequence
  $((C_j(k))_{j=1}^\infty)_{k=1}^\infty\in\PP(1/k)$ was
  arbitrary. Namely, if~\eqref{eq:Riemann-like} did not hold, it would
  be possible to construct a sequence
  $((C_j(k))_{j=1}^\infty)_{k=1}^\infty\in\PP(1/k)$ for
  which~\eqref{eq:Riemann-like-2} fails.
\end{proof}

\begin{lemma}
  \label{lemma:convolution-approximation}
  Consider $\varphi\in C_c(\R^2)$ and let $\psi$ be as in
  Lemma~\ref{lemma:Riemann-type}. Then
  \begin{equation*}
    \sum_{k\in\Z^2}\psi(x-kh)h^2\varphi(kh)\to\psi*\varphi(x)
    \text{ as } h\to 0,
  \end{equation*}
  uniformly in $x\in\R^2$.
\end{lemma}
\begin{proof}
  We can estimate
  \begin{equation*}
    \Big|\psi*\varphi(x) - \sum_{k\in\Z^2}\psi(x-kh)h^2\varphi(kh)\Big|
    \le A + B,
  \end{equation*}
  where
  \begin{align*}
    A &:=\|\varphi\|_\infty\sum_{k\in\Z^2}
        \int_{kh+[0,h)^2}\big|\psi(x-y)-\psi(x-kh)\big|\,dy,\text{ and }\\
    B &:= \|\psi\|_\infty\sum_{k\in\Z^2}
        \int_{kh+[0,h)^2}\big|\varphi(y)-\varphi(kh)\big|\,dy.
  \end{align*}

  The sum in $A$ can be bounded from the above by
  \begin{equation*}
    \begin{split}
      &\sum_{k\in\Z^2} h^2\,\sup\big\{|\psi(z)-\psi(z')|
      : {z,z'}\in x-kh-[0,h)^2\big\}\\
      &\qquad\le \sup\Big\{\sum_{j=1}^\infty \lambda_2(C_j)
      \sup_{{z,z'}\in C_j}|\psi(z)-\psi(z')| :
      (C_j)_{j=1}^\infty\in\mathcal{P}(\sqrt{2}h) \Big\}.
    \end{split}
  \end{equation*}
  By Lemma~\ref{lemma:Riemann-type} this tends to zero as
  $h\to\infty$.
  
  The number of nonzero terms in $B$ is bounded from the above by
  $C/h^2$, where $C>0$ is a constant independent of $h$. Consequently,
  $B$ can be estimated from the above by
  $C'\sup\{|\varphi(z)-\varphi(z')|:|z-z'|^2\le 2h^2\}$, which tends
  to zero as $h\to 0$ by the uniform continuity of $\varphi$.
\end{proof}

\begin{proof}[Proof of Theorem~\ref{thm:UAT}]
  It is enough to consider the case with a single output ($k_0=1$),
  for the general case follows from a componentwise construction of
  $\mathcal{N}$.
  
  For any parameters $(a,b)\in\C^{i_0}\times U$ such that
  \begin{equation}
    \label{eq:smalness-condition}
    a\cdot z+b\in U \text{ for every } z\in\Bbar,
  \end{equation}
  define a bounded function $f_{a,b}:\Bbar\to\C$ by setting
  $f_{a,b}(z):=\rho(a\cdot z+b)$. Then define
  \begin{equation}
    \label{eq:Sigma}
    \Ssigma
    := \overbar{\linspan}\{f_{a,b}:
    (a,b)\in\C^{i_0}\times U\text{ satisfies~\eqref{eq:smalness-condition}}\}
    \subset\mathcal{B}(\Bbar).
  \end{equation}
  Here $\mathcal{B}(\Bbar)$ is the complex algebra of bounded
  functions on $\Bbar$ equipped with the supremum norm, and the
  closure of the span is with respect to that norm. The theorem will
  be proved by showing that $\Ssigma$ includes the subset of
  continuous functions of $\B(\Bbar)$.
  
  Let $\varphi$ be a mollifier on $\R^2$ and define
  $\varphi_p(s):=p^2\varphi(ps)$ for $p=1,2,\ldots$

  Fix an integer $m\ge 0$ and find open sets $V$ and $W$ such that
  $\emptyset\neq V\subset\subset W\subset\subset U$ and that
  $\Delta^m\rho$ does not vanish identically in $V$. Let
  $\chi\in C_c(U)$ be such that $\chi\equiv 1$ on $W$. The convolution
  \begin{equation*}
    (\chi\rho)*\varphi_p(s) := \int_{\R^2}(\chi\rho)(s-y)\varphi_p(y)\,dy
  \end{equation*}
  is then defined everywhere, and
  $(\chi\rho)*\varphi_p|_{V}\to\rho|_{V}$ as $p\to\infty$ in the sense
  of distributions in $V$.  Consequently, there exists an index $p_0$
  such that $V-\supp\varphi_{p_0}\subset W$ and
  $\Delta^m(\chi\rho)*\varphi_{p_0}$ does not vanish identically in
  $V$. Define $\stilde:\C\to\C$ by
  $\stilde(s) := (\chi\rho)*\varphi_{p_0}(s)$. Then $\stilde$ is
  smooth everywhere (in the sense of real differentiability), and
  $\Delta^m\stilde$ does not vanish identically in $V$.

  Fix $b\in V$ and choose $\epsilon>0$ such that if $a\in\C^{i_0}$ and
  $|a|_{\C^{i_0}}<\epsilon$, then $a\cdot z+b\in V$ for every
  $z\in\Bbar$. Denote $\N_0:=\{0,1,2,\ldots\}$, and for any
  multiindices ${\alpha,\beta}\in\N_0^{i_0}$ define
  \begin{equation}
    \label{eq:F}
    F_{\alpha,\beta}(a,z)
    := z^\alpha\,\zbar^\beta\,
    (\partial^{|\alpha|}\dbar^{|\beta|}\stilde)(a\cdot z + b),
  \end{equation}
  where $z\in\Bbar$ and $|a|_{\C^{i_0}}<\epsilon$. Here
  $z^\alpha=z_1^{\alpha_1}z_2^{\alpha_2}\cdots z_{i_0}^{\alpha_{i_0}}$
  (and analogously for $\zbar$, where the bar denotes elementwise
  complex conjugation), and $\partial:=(\partial_x-i\partial_y)/2$ and
  $\dbar:=(\partial_x+i\partial_y)/2$ are the Wirtinger derivatives
  operating on the complex function $\stilde(x+iy)$.

  If $|\alpha|=|\beta|=0$, then
  \begin{equation}
    \label{eq:induction}
    F_{\alpha,\beta}(a,\cdot)\in\Ssigma
    \text{ for every } a\text{ with }|a|_{\C^{i_0}}<\epsilon.
  \end{equation}
  Namely, suppose $|a|_{\C^{i_0}}<\epsilon$ and let $h\in\R$ and
  $k\in\Z^2$ be such that $\varphi_{p_0}(kh)\neq 0$. Then
  $a\cdot z+b-kh\in W$ for every $z\in\Bbar$, so the parameters
  $(a, b-kh)$ satisfy~\eqref{eq:smalness-condition}, and
  $\chi(a\cdot z+b-kh)=1$. Consequently,
  \begin{equation*}
    h^2\sum_{k\in\Z^2} \varphi_{p_0}(kh)f_{a,b-kh}(z) =
    \sum_{k\in\Z^2} (\chi\rho)(a\cdot z+b-kh)h^2\varphi_{p_0}(kh) \to
    F_{0,0}(a,z)
  \end{equation*}
  as $h\to 0$, uniformly in $z\in\Bbar$, by
  Lemma~\ref{lemma:convolution-approximation}, and therefore
  $F_{0,0}(a,\cdot)\in\Ssigma$.

  Next we will use Wirtinger calculus similarly
  to~\cite[Lemma~{4.2}]{voigtlaender2020universal} to show
  that~\eqref{eq:induction} holds for every $\alpha$ and $\beta$. For
  a function of $a\in\C^{i_0}$, let us denote by $\partial_{a_i}$ and
  $\dbar_{a_i}$ the partial Wirtinger derivatives with respect to the
  variable $a_i\in\C$. Fix ${\alpha,\beta}\in\N_0^{i_0}$, denote
  $F:=F_{\alpha,\beta}$, and assume that~\eqref{eq:induction} holds
  for $F$. The directional derivative of $F$ in the $a$-variable along
  a direction $v\in\C^{i_0}$, denoted by $(\partial/\partial v)F$,
  exists, and a calculation shows that
  \begin{equation}
    \label{eq:convergence}
    \frac{F(a+hv,z)-F(a,z)}{h}\to \frac{\partial}{\partial v}F(a,z)\text{ as } h\to 0, 
  \end{equation}
  uniformly in $z\in\Bbar$. For fixed $a$ and small $h\neq 0$, by
  assumption the left-hand side of~\eqref{eq:convergence} as a
  function of $z$ is in $\Ssigma$. Because of the uniform convergence
  and closedness of $\Ssigma$, also the right-hand side
  of~\eqref{eq:convergence} is in $\Ssigma$. It follows that
  $\partial_{a_i} F(a,\cdot)\in\Ssigma$ and
  $\dbar_{a_i} F(a,\cdot)\in\Ssigma$, for every
  $i=1,2,\ldots,{i_0}$. But by the chain rule for the Wirtinger
  derivatives,
  \begin{align*}
    \partial_{a_i}F(a,z)
    &= z_iz^\alpha\zbar^\beta(\partial\partial^{|\alpha|}\dbar^{|\beta|}\stilde)(a\cdot z + b)
      = F_{\alpha+e_i,\beta}(a,z), \text{ and }\\
    \dbar_{a_i}\partial F(a,z)
    &= \zbar_iz^\alpha\zbar^\beta(\dbar\partial^{|\alpha|}\dbar^{|\beta|}\stilde)(a\cdot z + b)
      = F_{\alpha,\beta+e_i}(a,z).
  \end{align*}
  Consequently, \eqref{eq:induction} is true for every $\alpha$ and
  $\beta$.

  Because $\Delta^m\stilde=(4\partial\dbar)^m\stilde$ does not vanish
  identically in $V$, for every $\alpha$ and $\beta$ such that
  $|\alpha|\le m$ and $|\beta|\le m$ there exists
  $b_{\alpha,\beta}\in V$ such that
  $\partial^{|\alpha|}\dbar^{|\beta|}\stilde(b_{\alpha,\beta})\neq
  0$. Then \eqref{eq:F} and~\eqref{eq:induction} with $a=0$ and
  $b=b_{\alpha,\beta}$ imply that
  $z^\alpha\,\zbar^\beta\in\Ssigma$. Consequently, $\Ssigma$ contains
  all functions of the form
  \begin{equation}
    \label{eq:SW-polys}
    p(z)
    = \sum_{\substack{|\alpha|\le m,\\|\beta|\le m}}
    c_{\alpha\beta}z^\alpha\zbar^\beta,
  \end{equation}
  where $z\in\Bbar$, $m\in\N$ and $c_{\alpha\beta}\in\C$ are
  arbitrary. Functions of the form~\eqref{eq:SW-polys} form a
  self-adjoint algebra of continuous complex functions on the compact
  set $\Bbar$, and that algebra separates points on $\Bbar$ and
  vanishes at no point of $\Bbar$. By the Stone--Weierstrass theorem
  \cite{MR0385023} such an algebra contains all continuous complex
  functions in its uniform closure, and therefore so does $\Ssigma$.
\end{proof}


\bibliographystyle{siam}%
\bibliography{refs}

\end{document}